\documentclass[12pt,twoside]{article}    
\usepackage{amssymb}
\usepackage{amsmath}
\usepackage{amsfonts}
\usepackage{graphicx}

\usepackage{url}

\providecommand{\U}[1]{\protect\rule{.1in}{.1in}}

\newtheorem{theorem}{Theorem}[section]
\newtheorem{remark}[theorem]{Remark}
\newtheorem{lemma}[theorem]{Lemma}

\newtheorem{definition}[theorem]{Definition}
\newtheorem{proposition}[theorem]{Proposition}
\newtheorem{example}[theorem]{Example}
\numberwithin{equation}{section}
\newenvironment{proof}[1][Proof]{\textbf{#1.} }{\ \rule{0.5em}{0.5em}}

\makeatletter \@addtoreset{equation}{section} \makeatother
\begin{document}


\pagestyle{myheadings}
\markboth{\hfill {\small A. Al-Hussein, A. Ninouh, B. Gherbal} \hfill}{\hfill
{\small  McKean-Vlasov FBDSDEs and applications to stochastic control} \hfill}


\thispagestyle{plain}

\begin{center}
{\large \textbf{McKean-Vlasov forward-backward doubly
stochastic differential equations and applications to stochastic control}} \\
\vspace{0.7cm} {\large AbdulRahman Al-Hussein$^{a}$, Abdelhakim Ninouh$^{b}$, Boulakhras Gherbal}$^{c,}$\footnote{This work is supported by the Algerian PRFU, project No. C00L03UN070120220005.}
\\
\vspace{0.2cm} {\footnotesize
{\it $^{a}$Department of Mathematics, College of Science, Qassim University, \\
 P.O.Box 6644, Buraydah 51452, Saudi Arabia \\ {\emph E-mail:} alhusseinqu@hotmail.com, hsien@qu.edu.sa \\ \smallskip $^{b,c}$Laboratory of Mathematical Analysis, Probability and Optimization, \\ University of Mohamed Khider,  P.O.Box 145, Biskra 07000, Algeria
 \\ {\emph E-mail:} math.n.abdalhkim@gmail.com, b.gherbal@univ-biskra.dz, bgherbal@yahoo.fr}} \\ \smallskip
\end{center}

\vspace{0.15cm}

\begin{abstract}
This paper investigates first the existence and uniqueness of solutions for McKean-Vlasov forward-backward doubly stochastic differential equations (MV-FBDSDEs) in infinite-dimensional real separable Hilbert spaces. These equations combine the features of forward-backward doubly stochastic differential equations with the mean-field approach, allowing the coefficients to depend on the solution distribution. We establish the existence and uniqueness of solutions for MV-FBDSDEs using the method of continuation and provide an example and a counterexample to illustrate our findings. 
Moreover, we extend the practical applicability of our results by employing them within the context of the stochastic maximum principle for a control problem governed by MV-FBDSDEs. This study  contributes to the field of stochastic control problems and presents the first analysis of MV-FBDSDEs in infinite-dimensional spaces.
\end{abstract}

\textbf{MSC 2010:} 60H10, 93E20.

\textbf{Keywords:}
Continuation method, cost functional, existence, forward-backward doubly SDEs, Hilbert space, McKean-Vlasov, maximum principle, monotonicity condition, optimal control, sufficient conditions, uniqueness.

\section{Introduction}
Pardoux and Peng \cite{PP} introduced
backward doubly stochastic differential equation (BDSDE) in 1994 to give probabilistic interpretation for the solutions of a class of semilinear stochastic PDEs. Since then, the theory of BDSDEs has developed  and found applications in various fields, including stochastic control, stochastic PDEs, and finance.

Motivated by BDSDEs, there has been a growing interest in doubly stochastic optimal control problems (see e.g. \cite{BCY1, WL}). Stochastic Hamilton systems, derived from the stochastic maximum principle of stochastic optimal control problems, fall under the category of forward-backward doubly stochastic differential equations (FBDSDEs). 
The existence and uniqueness of solutions for these equations, which can be fully coupled, have been studied in various works such as \cite{A, AG2, PS}, along with references therein. Peng and Shi \cite{PS} established the existence and uniqueness of FBDSDE solutions under certain monotone assumptions using the method of time continuation. Zhu et al. \cite{ZSG} extended the results of \cite{PS} to FBDSDEs in different dimensional Euclidean spaces, relaxing the imposed monotonicity assumptions. Additionally, Al-Hussein and Gherbal \cite{AG2} studied FBDSDEs with Poisson jumps, while Al-Hussein \cite{A} explored FBDSDEs in infinite dimensions.

Mean-field stochastic differential equations (SDEs), also known as SDEs of McKean-Vlasov type, represent another type of SDEs where the coefficients can depend on the distribution of the solution, as shown in \cite{BYZ, CD2} and the references therein. In accordance to Lasry and Lions \cite{Las-Lio} and the related references therein, these equations have been widely used in finance, quantum chemistry, and game theory.
Mean-field backward stochastic differential equations, called also BSDEs of McKean-Vlasov type (MV-BSDEs), were introduced by Buckdahn et al. \cite{BDLP} as the mean square limit of an interacting particle system of BSDEs.

It is worth knowing that the stochastic maximum principle approaches to the solutions of optimal control problems for mean field SDEs naturally reduce to the solutions of mean field FBSDE systems; cf. e.g. \cite{BYZ, CD3, CD1, CD2, Song-Wang}. The existence of solutions for MV-BSDEs and McKean-Vlasov FBSDEs (MV-FBSDEs) has been investigated in various works, including \cite{BYZ, CD1, CD2}, along with relevant references therein. Additionally, the works \cite{Ahm} and \cite{Mah-Mc} provide insights into McKean-Vlasov equations in Hilbert spaces and their applications.

In this paper, we have two main objectives. Firstly,  we aim to establish the existence and uniqueness of the solution for the following McKean-Vlasov forward-backward doubly stochastic differential equations (MV-FBDSDEs):
\begin{equation}\label{eq:2.1}
\left \{
\begin{array}[c]{l}
dy_{t}=f\left(  t,y_{t},Y_{t},z_{t},Z_{t},\mathbb{P}_{\left(  y_{t}
,Y_{t},z_{t},Z_{t}\right)  }\right)  dt\\
\hspace{3.75cm} + \, g(  t,y_{t},Y_{t},z_{t},Z_{t},\mathbb{P}_{\left(y_{t},Y_{t},z_{t},Z_{t}\right)  })  \, dW_{t}-z_{t}\,\overleftarrow{d B}_{t}, \medskip \\
dY_{t}=F\left(  t,y_{t},Y_{t},z_{t},Z_{t},\mathbb{P}_{\left( y_{t},Y_{t},z_{t},Z_{t}\right)  }\right)  dt \\
\hspace{3.75cm} + \, G(  t,y_{t},Y_{t},z_{t},Z_{t},\mathbb{P}_{\left(y_{t},Y_{t},z_{t},Z_{t}\right)})  \, \overleftarrow{d B}_{t}+Z_{t}\,
dW_{t}, \medskip \\
y_{0}=x,\text{ }Y_{T}=h\left(  y_{T},\mathbb{P}_{y_{T}}\right) .
\end{array}
\right.
\end{equation}
We consider these equations in infinite dimensional real separable Hilbert spaces.
The system (\ref{eq:2.1}) incorporates mutually independent cylindrical Wiener processes $\left( W_{t}\right)_{t\geq0}$ and $\left(B_{t}\right)_{t\geq0}$  on real separable Hilbert spaces $E_{1}$ and $E_{2}$, respectively. The mappings $f, g, F, G$ are allowed to depend on all random variables $\left(y,Y,z,Z\right)$ in addition to their distribution $\mathbb{P}_{\left(y,Y,z,Z\right)}$, thereby enhancing the generality of the system, besides being fully-coupled.

Secondly, we demonstrate that this work contributes to laying a solid foundation for studying stochastic control problems governed by MV-FBDSDEs. Specifically, in
 Section~\ref{sec:appl}, we apply the results here to the stochastic maximum principle for MV-FBDSDEs. As is well-known, dynamic programming requires the solution to satisfy the Markov property, which does not hold in general due to the presence of distributions in the system~(\ref{eq:2.1}). Therefore, the maximum principle remains the suitable tool for studying such control problems. To the best of our knowledge, our present work is the first to address MV-FBDSDEs in infinite-dimensional spaces and their applications to stochastic optimal control.

The paper is organized as follows: Section~\ref{sec:2} introduces the problem formulation, presenting  MV-FBDSDEs and stating the assumptions on the coefficients.
In Section~\ref{sec:3}, we rigorously establish the existence and uniqueness of the solution for MV-FBDSDEs~(\ref{eq:2.1}), accompanied by compelling proofs. Additionally, at the end of Section~\ref{sec:3}, we provide an illustrative example and counterexample that shed light on the implications of our results. Finally, in Section~\ref{sec:4}, we demonstrate the practical applications of MV-FBDSDEs by applying them to stochastic optimal control.

\section{Notation and Formulation of the Problem}\label{sec:2}
Consider a complete probability space $\left(\Omega,\mathcal{F},\mathbb{P}\right)$ with a fixed time duration $T>0$.  The class of $\mathbb{P}$-null sets of $\mathcal{F}$ is denoted as $\mathcal{N}$. Let $E_{1}$ and $E_{2}$ be real and separable Hilbert spaces. We suppose that $\left(  W_{t}\right)  _{t\geq0}$ and 
$\left(B_{t}\right)_{t\geq0}$ are two mutually independent cylindrical Wiener
processes on $E_{1}$ and $E_{2}$, respectively. For each $t\in \left[ 0,T\right]$,
we define the $\sigma$-algebra $\mathcal{F}_{t}:=\mathcal{F}_{t,T}^{B}\vee \mathcal{F}_{t}^{W}$, which is generated by  $\mathcal{F}_{t,T}^{B}\cup \mathcal{F}_{t}^{W}$. Here, 
$\mathcal{F}_{s,t}^{\theta}=\sigma \left(\theta_{r}-\theta_{s},s\leq r\leq t\right)  \vee \mathcal{N}$ and $\mathcal{F}_{t}^{\theta}=\mathcal{F}_{0,t}^{\theta}$, for any process 
$\theta_{t}$. The collection $\left(\mathcal{F}_{t}\right)_{0\leq t\leq T}$ is neither increasing, nor decreasing, and thus does not form a filtration on 
$(\Omega, \mathcal{F})$.

We shall investigate systems governed by nonlinear MV-FBDSDEs. These systems are described by the equations, presented in (\ref{eq:2.1}).
In these equations, the integral with respect to $\overleftarrow{d B}_{t}$ represents  a backward It\^{o} integral, while the integral with respect to
$dW_{s}$ is a standard forward It\^{o} integral. These two types of integrals
are particular cases of It\^{o}-Sokorohod integral. Here, for  a random variable $X$ in a separable Hilbert space, 
$\mathbb{P}_{X}$ denotes the probability measure induced by $X$. The term ``nonlinear" used to describe the system (\ref{eq:2.1}) refers not only to the fact to the fact that the coefficients $f$, $g$, $F$, and $G$ could be nonlinear functions of the vector process $(y_{t},Y_{t},z_{t},Z_{t})$ at time $t$, but also to the fact that they depend on its distribution $\mathbb{P}_{( y_{t},Y_{t},z_{t},Z_{t})}$.

If $S$ is a separable real Hilbert space with norm $\| \cdot \|$, we denote by $\mathcal{P}\left(  S\right)$ to the space of all probability measures on $(S, \mathcal{B}(S))$, and by $\mathcal{P}_{2}\left(  S\right)$ to the subspace of
$\mathcal{P}\left(  S\right)$ of all probability measures having finite second order moments on $S$. We endow $\mathcal{P}_{2}\left(S\right)$ with the
\emph{$2$-Wasserstein distance} as follows:
\begin{align}\label{eq:2.2a}
\bar{w}_{2}\left(  \mu_{1},\mu_{2}\right)  & =  \inf \! \Bigg\{\!\Big(
\int \limits_{S\times S}
\|x-y\|^{2}  \lambda (dx,dy)  \Big)^{\frac{1}{2}}
| \, \lambda \in \mathcal{P} \left(  S\times S\right)  \text{with
marginals } \mu_{1},\mu_{2} \Bigg\} \nonumber \\
&=\inf \Big\{ \left(  \mathbb{E}\left[  \| X-Y\|^{2}\, \right]  \right)
^{\frac{1}{2}} \, | \; X,Y: \Omega
\rightarrow S   \text{ with
}\mathbb{P}_{X_{1}}=\mu_{1},\mathbb{P}_{X_{2}}=\mu_{2} \Big\} .
\end{align}
This definition makes $\mathcal{P}_{2}\left(S\right)$ a complete separable metric space.
We observe that if $X_{1}$ and $X_{2}$ are two square integrable random
variables taking their values in $S$, then the following inequality holds:
\begin{equation}\label{eq:2.2}
\| \, \mathbb{E}\left[  X_{1}\right]  -\mathbb{E}\left[  X_{2}\right]
\|  \leq \bar{w}_{2}\left(  \mathbb{P}_{X_{1}},\mathbb{P}_{X_{2}}\right)  \leq \left(\mathbb{E}\left[  \|X_{1}-X_{2}\|^{2}\, \right] \right)^{\frac{1}{2}} .
\end{equation}

Let $H$ be a separable real Hilbert space $H$ with inner product
$\left \langle \cdot,\cdot \right \rangle _{H}$ and norm
$\left \vert \cdot \right \vert _{H}$. We denote by $L_2 \left(  E_{i},H\right)$ to the space of all Hilbert-Schmidt operators from $E_{i}$ into $H$, where $i=1,2$. The inner product on $L_2\left( E_{i},H\right)$ is denoted by
$\left \langle \cdot ,\cdot \right \rangle _{L_2\left(E_{i},H\right)}$, and the norm induced by this inner product is denoted by
$\left \Vert \cdot \right \Vert_{L_2\left(  E_{i},H\right)}$.
For any $v^{1}=\left(  y^{1},Y^{1},z^{1},Z^{1}\right)$,
$v^{2}=\left(y^{2},Y^{2},z^{2},Z^{2}\right)  \in \mathbb{H}^2:=H\times H\times
L_2\left(  E_{2},H\right)  \times L_2\left(  E_{1},H\right)$, we define
$$\left(  v^{1},v^{2}\right)  =\left \langle y^{1},y^{2}\right \rangle
_{H}+\left \langle Y^{1},Y^{2}\right \rangle _{H}+\left \langle z^{1},z^{2}\right \rangle _{L_2\left(  E_{2},H\right)  }+\left \langle
Z^{1},Z^{2}\right \rangle _{L_2\left( E_{1},H\right)} ,$$
 and let
 $\left \vert v^{1}\right \vert
=\left[  \left(  v^{1},v^{1}\right) \right]^{\frac{1}{2}}$ be its norm.
Finally, for a separable Hilbert space $E$, we denote by $\mathfrak{M}^{2}\left(  \left[
0,T\right]  ,E\right)  $ to the space of all $E$-valued stochastic processes $\left(  X_{t}\right)_{0\leq t\leq T}$ such that for each $t\in [0,T]$, $X_t$ is $\mathcal{F}_{t}$-measurable, and $\mathbb{E}\, [\int_{0}^{T}|X_{t}|_{E}^{2}\, dt]<+\infty.$
Then it is evident that $\mathfrak{M}^{2}\left(  \left[  0,T\right]  ,E\right)$ is a Hilbert space endowed with the canonical norm
$$\left \Vert X\right \Vert
=\Big(\mathbb{E}\, [\int_{0}^{T}|X_{t}|_{E}^{2}\, dt]\Big)^{1/2} .$$

\begin{definition}
A quadruple  
$\left(  y,Y,z,Z\right)\in \mathfrak{M}^{2}\left(  \left[  0,T\right]  ,\mathbb{H}^2 \right)$ is called a \emph{solution} of MV-FBDSDEs~(\ref{eq:2.1}), if it satisfies 
($\mathbb{P}$-almost surely) the following integral systems for each $t\in [0, T]$:
\[
\left \{
\begin{array}[c]{l}
y_{t}=x+\int_{0}^{t}f\left(  s,y_{s},Y_{s},z_{s},Z_{s},\mathbb{P}_{\left(y_{s},Y_{s},z_{s},Z_{s}\right)  }\right)  ds\\
\hspace{3.5cm}+\int_{0}^{t}g\left(  s,y_{s},Y_{s},z_{s},Z_{s},\mathbb{P}_{\left(  y_{s},Y_{s},z_{s},Z_{s}\right)  }\right)  dW_{s}-\int_{0}^{t}z_{s}\, \overleftarrow{d B}_{s} , \medskip \\
Y_{t}=h\left(  y_{T},\mathbb{P}_{y_{T}}\right)  -\int_{t}^{T}F\left(
s,y_{s},Y_{s},z_{s},Z_{s},\mathbb{P}_{\left(  y_{s},Y_{s},z_{s},Z_{s}\right)
}\right)  ds\\
\hspace{3.5cm}-\int_{t}^{T}G\left(  s,y_{s},Y_{s},z_{s},Z_{s},\mathbb{P}_{\left(  y_{s},Y_{s},z_{s},Z_{s}\right)  }\right)   \overleftarrow{d B}_{s}-\int_{t}^{T}Z_{s}\, dW_{s}.
\end{array}
\right.
\]
\end{definition}

\medskip

For convenience, we introduce the notation:
\[
\begin{array}[c]{ll}
v=\left(  y,Y,z,Z\right), \;\, A\left(  t,v,\mu \right)  =\left(
F,f,G,g\right)  \left(  t,v,\mu \right),  & \\
\left(  A,v\right)  =\left \langle F,y\right \rangle _{H}+\left \langle
f,Y\right \rangle _{H}+\left \langle G,z\right \rangle _{L_2\left(E_{2},H\right)}+\left \langle g,Z\right \rangle _{L_2\left(E_{1},H\right)} ,
\end{array}\]
where $\mu$ is a probability measure on $\mathbb{H}^2$.

We now state our main assumptions for the mappings: 
\begin{center}
$\left(  f,F\right)  :\Omega\times [0, T]\times\mathbb{H}^2 \times \mathcal{P}_{2}\left( \mathbb{H}^2 \right)  \rightarrow H$, \\ $g:\Omega\times [0, T]\times\mathbb{H}^2 \times \mathcal{P}_{2}\left( \mathbb{H}^2 \right)
\rightarrow L_2\left(  E_{1},H\right)$, \\
$G:\Omega\times [0, T]\times\mathbb{H}^2 \times \mathcal{P}_{2}\left( \mathbb{H}^2 \right)
\rightarrow L_2\left(  E_{2},H\right)$, \\ $h:\Omega\times H\times \mathcal{P}_{2}\left( \mathbb{H}^2 \right) \rightarrow H$.
\end{center}

\medskip

$\left(\mathbf{A}_{1}\right)$ (\textbf{Lipschitz conditions}) There exist  $C>0, \, 0<\gamma<1/2$ such that for each
$\left(  v^{i},\mu^{i}\right)  :=\left(  y^{i},Y^{i},z^{i},Z^{i},\mu^{i}\right)    \in \mathbb{H}^2 \times \mathcal{P}_{2}\left(
\mathbb{H}^2 \right)$, if denoting $v^{i,z}=\left(
y^{i},Y^{i},Z^{i}\right)$ and  $v^{i,Z}=\left(  y^{i},Y^{i},z^{i}\right)$ as elements of $H\times H\times L_2\left(
E_{1},H\right)$ and $H\times
H\times L_2\left(  E_{2},H\right)$, respectively for $i=1,2,$ then we have for each
$t\in \left[  0,T\right]$,
\[%
\begin{array}
[c]{ll}%
{\rm\left(  i\right)}  \left \vert \left( f,F\right)  \left(  t,v^{1},\mu^{1}
\right)  -\left( f,F\right)  \left(  t,v^{2},\mu^{2}\right)  \right \vert
_{H}\leq C\left(  \left \vert v^{1}-v^{2}\right \vert +\bar{w}_{2}\left(
\mu^{1},\mu^{2}\right)  \right)  & \smallskip \\
{\rm\left(  ii\right)} \left \Vert G\left(  t,v^{1},\mu^{1}\right)  -G\left(
t,v^{2},\mu^{2}\right)  \right \Vert _{L_2\left(  E_{2},H\right)
}^{2}\leq C \left \vert v^{1,Z}-v^{2,Z}\right \vert ^{2}  & \\
\hspace{7.5cm}+\, \gamma \left(  \left \Vert Z^{1}-Z^{2}\right \Vert _{L_2
\left(  E_{1},H\right)  }^{2}+\bar{w}_{2}^2\left(  \mu^{1},\mu^{2}\right)\right)  & \smallskip \\
{\rm\left(  iii\right)} \left \Vert g\left(  t,v^{1},\mu^{1}\right)  -g\left(
t,v^{2},\mu^{2}\right)  \right \Vert _{L_2\left(  E_{1},H\right)
}^{2}\leq C  \left \vert v^{1,z}-v^{2,z}\right \vert ^{2}  & \\
\hspace{7.5cm}+\, \gamma \left(  \left \Vert z^{1}-z^{2}\right \Vert _{L_2
\left(  E_{2},H\right)  }^{2}+\bar{w}_{2}^2\left(  \mu^{1},\mu^{2}\right) \right)  & \smallskip \\
{\rm\left(  iv\right)} \left \vert h\left(  y^{1},\mu^{1}\right)  -h\left(
y^{2},\mu^{2}\right)  \right \vert _{H}\leq C \left(  \left \vert y^{1}
-y^{2}\right \vert _{H}+\bar{w}_{2} \left(  \mu^{1},\mu^{2}\right)
\right)  . &
\end{array}
\]

$\left(\mathbf{A}_{2}\right)$ (\textbf{Monotonicity conditions}) Assume that there exist non-negative constants $\theta_{1},$
$\theta_{2}$, and $\alpha_{1}$ with $\theta_{1}+\theta_{2}>0,$ $\alpha
_{1}+\theta_{2}>0$ such that for any random variables $v^{1} :=\left(  y^{1},Y^{1},z^{1},Z^{1}\right)$ and $v^{2} := \left(y^{2},Y^{2},z^{2},Z^{2}\right)$ taking values in $\mathbb{H}^{2}$ and for any $t\in \left[  0,T\right]$,
we have
\begin{align*}
& \hspace{-0.75cm}{\rm\left(i\right)} \;  \, \mathbb{E} \left[ \Big(  A\left(  t,\upsilon^1,\mathbb{P}_{\left(y^{1},Y^{1},z^{1},Z^{1}\right)}\right)  -A\left(
t,\upsilon^2 , \mathbb{P}_{\left(y^{2},Y^{2},z^{2},Z^{2}\right)} \right)  ,\upsilon^1 - \upsilon^2 \Big)  \right]  \\
& \hspace{2cm}  \leq
 -\, \theta_{1}\, \mathbb{E} \left[ \left \vert y^{1}-y^{2}\right \vert_{H}^{2}+\left \Vert z^{1}-z^{2}\right \Vert _{L_2\left(  E_{2},H\right)}^{2} \right]
\\
& \hspace{3.5cm}
 - \, \theta_{2}\, \mathbb{E} \left[  \left \vert Y^{1}-Y^{2}\right \vert
_{H}^{2}+\left \Vert Z^{1}-Z^{2}\right \Vert _{L_2\left(E_{1},H\right)  }^{2}\right] ,
 \\
 &  \hspace{-0.75cm} {\rm\left(ii\right)} \; \, \mathbb{E} \left[  \left \langle h\left(  y^{1},\mathbb{P}_{y^1}\right)
  -h\left(  y^{2},\mathbb{P}_{y^2}\right)  ,y^{1}-y^{2}\right \rangle_{H} \right] \geq \alpha_{1} \, \mathbb{E} \left[  \left \vert y^{1}-y^{2}\right \vert_{H}^{2} \right].
\end{align*}

$\left(\mathbf{A}_{3}\right)$
For each element $v=(y,Y,z,Z)$ of $\mathbb{H}^2 $ and for each $\mu\in \mathcal{P}_{2}\left(  \mathbb{H}^2 \right)$, we have ${A\left(
\cdot,v,\mu \right)  \in \mathfrak{M}^{2}\left(  \left[  0,T\right]
,\mathbb{H}^2 \right) }$ and $h\left(  y,\mu \right)  \in L^{2}\left(  \Omega
,\mathcal{F}_{T},\mathbb{P},H\right)$.

\bigskip

The observation that Wasserstein's distance of two probability measures
is bounded below by the Euclidean norm of the difference of their respective
expectations, as demonstrated in (\ref{eq:2.2}), motivates considering the same
research problem under different influences of various Lipschitz constraints.

\begin{remark}\label{rk:rk-on-assumptions}
{\rm (i)} As a special case, when $h$ does not depend on $\left( y,\mu \right)$,
i.e. ${h\left( y,\mu \right) =\xi}$ for a given
$\xi \in L^{2}\left(\Omega ,\mathcal{F}_{T},\mathbb{P},H\right)$,
the two monotonicity conditions imposed in $\left(\mathbf{A}_{2}\right)$ collapse to the following condition:
\begin{align*}
& \hspace{-0.5cm} \mathbb{E} \left[ \Big(  A\left(  t,\upsilon^1,\mathbb{P}_{\left(y^{1},Y^{1},z^{1},Z^{1}\right)}\right)  -A\left(
t,\upsilon^2 , \mathbb{P}_{\left(y^{2},Y^{2},z^{2},Z^{2}\right)} \right)  ,\upsilon^1 - \upsilon^2 \Big)  \right]  \\
& \hspace{2cm}  \leq
 -\, \theta_{1}\, \mathbb{E} \left[ \left \vert y^{1}-y^{2}\right \vert_{H}^{2}+\left \Vert z^{1}-z^{2}\right \Vert _{L_2\left(  E_{2},H\right)}^{2} \right]
\\
& \hspace{4cm}
 - \, \theta_{2}\, \mathbb{E} \left[  \left \vert Y^{1}-Y^{2}\right \vert
_{H}^{2}+\left \Vert Z^{1}-Z^{2}\right \Vert _{L_2\left(E_{1},H\right)  }^{2}\right]
\end{align*}
for some constants $\theta_{1} \geq 0$ and $\theta_{2} > 0$.
\smallskip \\
{\rm (ii)} The assumption $\left(\mathbf{A}_{2}\right)$ can be replaced by the following ones, while maintaining essentially the same proofs for the theorems in the following section and their respective lemmas.

$\left(\mathbf{A}_{2}\right)':$\; $\forall\; v^{1} :=\left(  y^{1},Y^{1},z^{1},Z^{1}\right), v^{2} := \left(y^{2},Y^{2},z^{2},Z^{2}\right)\in \mathbb{H}^{2}$ and $\forall\; t\in \left[  0,T\right]$,
\begin{align*}
& \hspace{-0.75cm} \, \mathbb{E} \left[ \Big(  A\left(  t,\upsilon^1,\mathbb{P}_{\left(y^{1},Y^{1},z^{1},Z^{1}\right)}\right)  -A\left(
t,\upsilon^2 , \mathbb{P}_{\left(y^{2},Y^{2},z^{2},Z^{2}\right)} \right)  ,\upsilon^1 - \upsilon^2 \Big)  \right]  \\
& \hspace{2cm}  \geq
\theta_{1}\, \mathbb{E} \left[ \left \vert y^{1}-y^{2}\right \vert_{H}^{2}+\left \Vert z^{1}-z^{2}\right \Vert _{L_2\left(  E_{2},H\right)}^{2} \right]
\\
& \hspace{3.5cm}
 + \, \theta_{2}\, \mathbb{E} \left[  \left \vert Y^{1}-Y^{2}\right \vert
_{H}^{2}+\left \Vert Z^{1}-Z^{2}\right \Vert _{L_2\left(E_{1},H\right)  }^{2}\right] 
\end{align*}
and 
\begin{align*}
 &  \mathbb{E} \left[  \left \langle h\left(  y^{1},\mathbb{P}_{y^1}\right)
  -h\left(  y^{2},\mathbb{P}_{y^2}\right)  ,y^{1}-y^{2}\right \rangle_{H} \right] \leq - \, \alpha_{1} \, \mathbb{E} \left[  \left \vert y^{1}-y^{2}\right \vert_{H}^{2} \right].
\end{align*}
\end{remark}

\section{Existence and Uniqueness Theorems}\label{sec:3}
In this section, we establish our main result of the existence
and uniqueness of the solution to MV-FBDSDEs, which is a system of nonlinear fully coupled FBDSDEs of McKean-Vlasov type.

\subsection{Uniqueness of the Solutions of MV-FBDSDEs~(\ref{eq:2.1})}\label{subsec:uniq-sec}
The following theorem gives conditions that guarantee the uniqueness of the solution of MV-FBDSDEs~(\ref{eq:2.1}).
\begin{theorem}
\label{thm:main thm1} Under
 $\left(\mathbf{A}_{1}\right)$--$\left(\mathbf{A}_{3}\right)$, MV-FBDSDEs~(\ref{eq:2.1}) has at most one solution $\left(y,Y,z,Z\right)$ in
  $\mathfrak{M}^{2}\left( \left[ 0,T\right] ,\mathbb{H}^2\right)$.
\end{theorem}

Let us begin by introducing the integration by parts formula, commonly referred to as 
It\^{o}'s formula. This formula is derived from the classical It\^{o}'s formula, as it can be gleaned, for instance, from \cite{PP}.
\begin{proposition}\label{prop:Ito formula}
Let $\left( \alpha,\beta,\gamma,\delta \right)$
and
$\left(\tilde{\alpha},\tilde{\beta},\tilde{\gamma},\tilde{\delta}\right)$ be elements of
 $\mathfrak{M}^{2}\left( \left[ 0,T\right] ,\mathbb{H}^2 \right)$
and $\mathfrak{M}^{2}\left( \left[ 0,T\right] ,\mathbb{H}^2 \right)$, respectively. Assume that
\begin{equation*}
\left \{
\begin{array}{ll}
\alpha_{t}=\alpha_{0}+\int_{0}^{t}\beta_{s}\, ds+\int_{0}^{t}\delta _{s}\,
\overleftarrow{d B}_{s}+\int_{0}^{t}\gamma_{s}\, dW_{s} , &  \smallskip \\
\tilde{\alpha}_{t}=\tilde{\alpha}_{0}+\int_{0}^{t}\tilde{\beta}
_{s}\, ds+\int_{0}^{t} \tilde{\delta}_{s}\, \overleftarrow{d B}_{s}+\int_{0}^{t}
\tilde{\gamma}_{s}\, dW_{s} , &
\end{array}
\right.
\end{equation*}
for all $t\in \left[ 0,T\right]$. Then, for each $t\in \left[ 0,T\right]$,
\begin{equation*}
\left \langle \alpha_{t},\tilde{\alpha}_{t}\right \rangle _{H}=\left \langle
\alpha_{0},\tilde{\alpha}_{0}\right \rangle _{H}+\int_{0}^{t}\left \langle
\alpha_{s},d \tilde{\alpha}_{s}\right \rangle _{H}+\int_{0}^{t}\left \langle
\tilde{\alpha}_{s},d\alpha_{s}\right \rangle _{H}+\int_{0}^{t}d\left \langle
\alpha_{s},\tilde{\alpha}_{s}\right \rangle _{H} \quad \mathbb{P}-a.s.,
\end{equation*}
and
\begin{equation*}
\begin{array}{ll}
\mathbb{E}\left[ \left \langle \alpha_{t},\tilde{\alpha}_{t}\right \rangle
_{H}\right] =\mathbb{E}\left[ \left \langle \alpha_{0},\tilde{\alpha}
_{0}\right \rangle _{H}\right] +\mathbb{E}\left[ \int_{0}^{t}\left \langle
\alpha_{s},d\tilde{\alpha}_{s}\right \rangle _{H}\right] +\mathbb{E}\left[
\int_{0}^{t}\left \langle \tilde{\alpha}_{s},d\alpha_{s}\right \rangle_{H}
\right] &  \\
\hspace{4.5cm} - \, \mathbb{E}\left[ \int_{0}^{t} \langle\, \delta_{s},
\tilde{\delta}_{s}\, \rangle_{L_2(E_2,H)} \, ds\right] +\mathbb{E}\left[ \int_{0}^{t}\left
\langle \gamma_{s},\tilde{\gamma}_{s}\right \rangle_{L_2(E_1,H)} \, ds\right] . &
\end{array}%
\end{equation*}
\end{proposition}

\medskip

\begin{proof}[Proof of Theorem~\protect\ref{thm:main thm1}]
Let $v^{i}=\left( y^{i},Y^{i},z^{i},Z^{i}\right)$, for $i=1,2$, be two
solutions of system~(\ref{eq:2.1}). To simplify the notation, we denote
\begin{equation*}
\begin{array}{ll}
\bigtriangleup v=\left( \bigtriangleup y,\bigtriangleup Y,\bigtriangleup
z,\bigtriangleup Z\right) =\left(
y^{1}-y^{2},Y^{1}-Y^{2},z^{1}-z^{2},Z^{1}-Z^{2}\right), &  \\
\bigtriangleup f_{t}=f\left( t,y_{t}^{1},Y_{t}^{1},z_{t}^{1},Z_{t}^{1},%
\mathbb{P}_{\left( y_{t}^{1},Y_{t}^{1},z_{t}^{1},Z_{t}^{1}\right) }\right)
-f\left( t,y_{t}^{2},Y_{t}^{2},z_{t}^{2},Z_{t}^{2},\mathbb{P}_{\left(
y_{t}^{2},Y_{t}^{2},z_{t}^{2},Z_{t}^{2}\right) }\right) , &  \\
\bigtriangleup g_{t}=g\left( t,y_{t}^{1},Y_{t}^{1},z_{t}^{1},Z_{t}^{1},%
\mathbb{P}_{\left( y_{t}^{1},Y_{t}^{1},z_{t}^{1},Z_{t}^{1}\right) }\right)
-g\left( t,y_{t}^{2},Y_{t}^{2},z_{t}^{2},Z_{t}^{2},\mathbb{P}_{\left(
y_{t}^{2},Y_{t}^{2},z_{t}^{2},Z_{t}^{2}\right) }\right), &  \\
\bigtriangleup F_{t}=F\left( t,y_{t}^{1},Y_{t}^{1},z_{t}^{1},Z_{t}^{1},%
\mathbb{P}_{\left( y_{t}^{1},Y_{t}^{1},z_{t}^{1},Z_{t}^{1}\right) }\right)
-F\left( t,y_{t}^{2},Y_{t}^{2},z_{t}^{2},Z_{t}^{2},\mathbb{P}_{\left(
y_{t}^{2},Y_{t}^{2},z_{t}^{2},Z_{t}^{2}\right) }\right), &  \\
\bigtriangleup G_{t}=G\left( t,y_{t}^{1},Y_{t}^{1},z_{t}^{1},Z_{t}^{1},%
\mathbb{P}_{\left( y_{t}^{1},Y_{t}^{1},z_{t}^{1},Z_{t}^{1}\right) }\right)
-G\left( t,y_{t}^{2},Y_{t}^{2},z_{t}^{2},Z_{t}^{2},\mathbb{P}_{\left(
y_{t}^{2},Y_{t}^{2},z_{t}^{2},Z_{t}^{2}\right) }\right) , &  \\
\bigtriangleup h_{T}=h\left( y_{T}^{1},\mathbb{P}_{y_{T}^{1}}\right)
-h\left( y_{T}^{2},\mathbb{P}_{y_{T}^{2}}\right) , &
\end{array}%
\end{equation*}
where $0\leq t\leq T.$

By applying It\^{o}'s formula (see Proposition~\ref{prop:Ito formula}) and
 $\left(\mathbf{A}_{2}\right) \mathrm{%
\left(i\right)}$ to $\left \langle
\bigtriangleup y_{t},\bigtriangleup Y_{t}\right \rangle _{H}$, we obtain
\begin{align*}
& \mathbb{E}\left[ \left \langle \bigtriangleup y_{T},\bigtriangleup
h_{T}\right \rangle _{H}\right] \\ &
\hspace{2cm} =\mathbb{E}\left[ \int_{0}^{T}\left( A\big( t,v_{t}^{1},
\mathbb{P}_{\left( y_{t}^{1},Y_{t}^{1},z_{t}^{1},Z_{t}^{1}\right) }\big)
-A\big( t,v_{t}^{2},\mathbb{P}_{\left(
y_{t}^{2},Y_{t}^{2},z_{t}^{2},Z_{t}^{2}\right) }\big) ,\bigtriangleup
v_{t}\right) dt\right] \\ &
\hspace{2cm}\leq
-\theta_{1}\, \mathbb{E}\left[ \int_{0}^{T}\big{(}\left \vert
y_{t}^{1}-y_{t}^{2}\right \vert _{H}^{2}+\left \Vert
z_{t}^{1}-z_{t}^{2}\right \Vert _{L_2\left( E_{2},H\right) }^{2}\big{)}\, dt\right]
\\ & \hspace{5cm} - \, \theta_{2}\, \mathbb{E}\left[ \int_{0}^{T} \big{(} \left \vert
Y_{t}^{1}-Y_{t}^{2}\right \vert _{H}^{2}+\left \Vert
Z_{t}^{1}-Z_{t}^{2}\right \Vert _{L_2\left( E_{1},H\right) }^{2}\big{)}\,  dt\right] .
\end{align*}
Hence, according to $\left(\mathbf{A}_{2}\right) \mathrm{\left(ii\right)}$, it follows that
\begin{eqnarray*}
& 0\leq \alpha_{1}\, \mathbb{E}\left[ \left \vert \bigtriangleup y_{T}\right
\vert _{H}^{2}\right] \leq
- \, \theta_{1}\, \mathbb{E}\left[ \int_{0}^{T}\big{(}\left \vert
y_{t}^{1}-y_{t}^{2}\right \vert _{H}^{2}+\left \Vert
z_{t}^{1}-z_{t}^{2}\right \Vert _{L_2\left( E_{2},H\right) }^{2}\big{)}\, dt\right]
\\ & \hspace{5cm} - \, \theta_{2}\, \mathbb{E}\left[ \int_{0}^{T} \big{(} \left \vert
Y_{t}^{1}-Y_{t}^{2}\right \vert _{H}^{2}+\left \Vert
Z_{t}^{1}-Z_{t}^{2}\right \Vert _{L_2\left( E_{1},H\right) }^{2}\big{)}\,  dt\right] .
\end{eqnarray*}

If both $\theta_{1}>0$ and  $\theta_{2}>0$ (e.g., when $\theta_{1}=\theta_{2}$), this inequality directly proves the uniqueness of $(y,Y,z,Z)$, so that we would not need to assume that $0<\gamma< 1/2$  for the purpose of establishing the uniqueness of the solutions of MV-FBDSDE~(\ref{eq:2.1}).
Therefore, let us consider the general case in the remaining part of the proof.

If $\theta_{2}>0$, we obtain ${\mathbb{E}\big[ \left \vert
Y_{t}^{1}-Y_{t}^{2}\right \vert _{H}^{2}\big] =0}$ and $\mathbb{E}\big[
\left \Vert Z_{t}^{1}-Z_{t}^{2}\right \Vert _{L_2\left( E_{1},H\right) }^{2}
\big] =0,$ which imply that $Y_{t}^{1}=Y_{t}^{2}$ and $Z_{t}^{1}=Z_{t}^{2} \; a.s.$
for all $0\leq t\leq T.$ Hence, according to (\ref{eq:2.1}), we have
\begin{equation*}
\bigtriangleup y_{t}=\int_{0}^{t}\bigtriangleup
\widehat{f}_{s}\, ds+\int_{0}^{t}\bigtriangleup \widehat{g}_{s}\, dW_{s}-\int_{0}^{t}\bigtriangleup
z_{s}\, \overleftarrow{d B}_{s} , \quad t\in \left[ 0,T\right],
\end{equation*}
where $\bigtriangleup \widehat{f}_{t}= \bigtriangleup f_{t}$ and $\bigtriangleup \widehat{g}_{t}= \bigtriangleup g_{t}$ in this case, or specifically
\begin{align*}
& \bigtriangleup \widehat{f}_{t}=f\left( t,y_{t}^{1},Y_{t}^{1},z_{t}^{1},Z_{t}^{1},
\mathbb{P}_{\left( y_{t}^{1},Y_{t}^{1},z_{t}^{1},Z_{t}^{1}\right) }\right)
-f\left( t,y_{t}^{2},Y_{t}^{1},z_{t}^{2},Z_{t}^{1},\mathbb{P}_{\left(
y_{t}^{2},Y_{t}^{1},z_{t}^{2},Z_{t}^{1}\right) }\right),  \\ &
\bigtriangleup \widehat{g}_{t}=g\left( t,y_{t}^{1},Y_{t}^{1},z_{t}^{1},Z_{t}^{1},
\mathbb{P}_{\left( y_{t}^{1},Y_{t}^{1},z_{t}^{1},Z_{t}^{1}\right) }\right)
-g\left( t,y_{t}^{2},Y_{t}^{1},z_{t}^{2},Z_{t}^{1},\mathbb{P}_{\left(
y_{t}^{2},Y_{t}^{1},z_{t}^{2},Z_{t}^{1}\right) }\right) .
\end{align*}

Applying It\^{o}'s formula  to $\left
\vert \bigtriangleup y_{t}\right \vert _{H}^{2}$ yields
\begin{align*}
\mathbb{E}\left[ \left \vert \bigtriangleup y_{t}\right \vert _{H}^{2}\right]
& +\mathbb{E}\big[ \int_{0}^{t} \left \Vert \bigtriangleup z_{s}\right \Vert
_{L_2\left( E_{2},H\right) }^{2}ds\big] &  \\
& =2\, \mathbb{E}\big[ \int_{0}^{t} \langle  \, \vert \!
\bigtriangleup \! \widehat{f}_{s} \vert _{H},\left \vert \bigtriangleup y_{s}\right
\vert _{H}\, \rangle _{H}ds\big] +\mathbb{E}\big[ \int_{0}^{t}\left
\Vert \bigtriangleup \widehat{g}_{s}\right \Vert _{L_2\left( E_{1},H\right) }^{2}ds%
\big]  \\
& \leq 2 \, \mathbb{E}\big[ \int_{0}^{t} \vert\! \bigtriangleup \!
\widehat{f}_{s} \vert _{H}\left \vert \bigtriangleup y_{s}\right \vert _{H}ds%
\big] +\mathbb{E}\big[ \int_{0}^{t}\left \Vert \bigtriangleup \widehat{g}_{s}\right
\Vert _{L_2\left( E_{1},H\right) }^{2}ds\big] .
\end{align*}
Hence, based on $\left( \mathbf{A}_{1}\right)$ and the inequality $ab\leq \frac{1}{2 \varepsilon}\, a^{2}+\frac{\varepsilon}{2}\, b^{2}$,
for any $\varepsilon > 0$, it follows that
\begin{align*}
\mathbb{E} \, \big[ \left \vert \bigtriangleup y_{t}\right \vert _{H}^{2}\big]
+\mathbb{E} \, \Big[ \int_{0}^{t}\left \Vert \bigtriangleup z_{s}\right \Vert
_{L_2\left( E_{2},H\right) }^{2}ds\Big] \leq 2 C \int_{0}^{t}\mathbb{E}\left[
\left( \left \vert \bigtriangleup y_{s}\right \vert _{H}^{2}+\left \Vert
\bigtriangleup z_{s}\right \Vert _{L_2\left( E_{2},H\right) }\left \vert
\bigtriangleup y_{s}\right \vert _{H}\right. \right. &  \\
\hspace{3cm}+\left. \left. \left \vert \bigtriangleup y_{s}\right \vert _{H}%
\big{(}\mathbb{E}\left[ \left \vert \bigtriangleup y_{s}\right \vert _{H}^{2}\right]\big{)}
^{\frac{1}{2}}+\left \vert \bigtriangleup y_{s}\right \vert _{H} \big{(} \mathbb{E}\,
\big[ \left \Vert \bigtriangleup z_{s}\right \Vert _{L_2\left(
E_{2},H\right) }^{2}\big{]} \big{)} ^{\frac{1}{2}}\right) ds \right] &  \\
\hspace{3cm}+\int_{0}^{t}\left(\left(C+\gamma\right)\, \mathbb{E}\left[ \left \vert \bigtriangleup
y_{s}\right \vert _{H}^{2}\right] +2\, \gamma \, \mathbb{E}\, \big[ \left \Vert
\bigtriangleup z_{s}\right \Vert _{L_2\left( E_{2},H\right) }^{2}\big]\right) ds
&  \\
\hspace{1.5cm}\leq 2 C \int_{0}^{t}\left( \left(2+\frac{1}{2\, \varepsilon}+\frac {1}{%
2\, \varepsilon}\right) \mathbb{E}\left[ \left \vert \bigtriangleup y_{s}\right
\vert _{H}^{2}\right] +\left( \frac{\varepsilon}{2}+\frac{\varepsilon }{2}\right)
\mathbb{E}\big[ \left \Vert \bigtriangleup z_{s}\right \Vert _{L_2\left(
E_{2},H\right) }^{2}\big]\right) ds &  \\
\hspace{3cm}+\int_{0}^{t}\left(\left(C+\gamma\right)\, \mathbb{E}\left[ \left \vert \bigtriangleup
y_{s}\right \vert _{H}^{2}\right] +2\gamma \, \mathbb{E}\, \big[ \left \Vert
\bigtriangleup z_{s}\right \Vert _{L_2\left( E_{2},H\right) }^{2}\big]\right) ds
&  \\
\hspace{1.5cm}\leq \int_{0}^{t}\left( \left(5C+\gamma+\frac{2C}{\varepsilon}\right) \mathbb{E}%
\left[ \left \vert \bigtriangleup y_{s}\right \vert _{H}^{2}\right] +\left(
2 \gamma+2 C\, \varepsilon \right) \mathbb{E}\, \big[ \left \Vert \bigtriangleup
z_{s}\right \Vert _{L_2\left( E_{2},H\right) }^{2}\big] \right)ds . &
\end{align*}
Since $0< \gamma< 1/2$, if we choose in this case $\varepsilon=\frac{1-2\gamma}{4C}$, we get
\begin{equation*}
\begin{array}{ll}
\hspace{-0.16cm} \mathbb{E} \, \big[ \left \vert \bigtriangleup y_{t}\right \vert _{H}^{2}\big]
+\left( \frac{1-2\gamma}{2}\right) \int_{0}^{t}\mathbb{E}\, \big[ \left \Vert
\bigtriangleup z_{s}\right \Vert _{L_2\left( E_{2},H\right) }^{2}\big]
ds\leq \left( 5C+\gamma+\frac{8C^{2}}{1-2\gamma}\right) \int_{0}^{t}\mathbb{E}\, \big[
\left \vert \bigtriangleup y_{s}\right \vert _{H}^{2}\big] ds. &
\end{array}
\end{equation*}

Now, Gronwall's inequality implies $y_{t}^{1}=y_{t}^{2} \;\, a.s.$ for all
$0\leq t\leq T$.  This, in turn, leads to
$$\int_{0}^{t}\mathbb{E}\, \big[ \left \Vert
\bigtriangleup z_{s}\right \Vert _{L_2\left( E_{2},H\right) }^{2}\big]
ds=0 ,$$ which yields $z_{t}^{1}=z_{t}^{2} \;\, a.s.$ for all
$0\leq t\leq T$.

If $\alpha _{1}>0$ and $\theta _{1}>0$, then we have
$$\mathbb{E}\big[ \left\vert
y_{t}^{1}-y_{t}^{2}\right\vert _{H}^{2}\big] =0, \;\; \mathbb{E}\left[
\left\vert \bigtriangleup y_{T}\right\vert _{H}^{2}\right]=0, \;\;\text{and}\;\;
\mathbb{E}\big[ \left\Vert z_{t}^{1}-z_{t}^{2}\right\Vert _{L_{2}\left(
E_{2},H\right) }^{2}\big] =0 .$$
As a result, $y_{t}^{1}=y_{t}^{2}$ and
$z_{t}^{1}=z_{t}^{2} \;\, a.s.$ for all $0\leq t\leq T$. Therefore,
 $$h\big(y_{T}^{1},P_{y_{T}^{1}}\big) =h\big( y_{T}^{2},P_{y_{T}^{2}}\big) ,$$ and so
\begin{equation*}
\bigtriangleup Y_{t}=-\int_{t}^{T}\bigtriangleup
\widehat{F}_{s}\, ds-\int_{t}^{T}\bigtriangleup \widehat{G}_{s}\,
\overleftarrow{d B}_{s}-\int_{t}^{T}\bigtriangleup Z_{s}\, dW_{s} ,
\end{equation*}
where here
\begin{align*}
\bigtriangleup \widehat{F}_{t}=F\left( t,y_{t}^{1},Y_{t}^{1},z_{t}^{1},Z_{t}^{1},%
\mathbb{P}_{\left( y_{t}^{1},Y_{t}^{1},z_{t}^{1},Z_{t}^{1}\right) }\right)
-F\left( t,y_{t}^{1},Y_{t}^{2},z_{t}^{1},Z_{t}^{2},\mathbb{P}_{\left(
y_{t}^{1},Y_{t}^{2},z_{t}^{1},Z_{t}^{2}\right) }\right) , &  \\
\bigtriangleup \widehat{G}_{t}=G\left( t,y_{t}^{1},Y_{t}^{1},z_{t}^{1},Z_{t}^{1},%
\mathbb{P}_{\left( y_{t}^{1},Y_{t}^{1},z_{t}^{1},Z_{t}^{1}\right) }\right)
-G\left( t,y_{t}^{1},Y_{t}^{2},z_{t}^{1},Z_{t}^{2},\mathbb{P}_{\left(
y_{t}^{1},Y_{t}^{2},z_{t}^{1},Z_{t}^{2}\right) }\right)  . &
\end{align*}
Next, apply It\^{o}'s formula to $\left\vert \bigtriangleup Y_{t}\right\vert
_{H}^{2}$ and utilize $\left( \mathbf{A}_{1}\right)$ to find that
\begin{align}\label{eq:3.1}
& \mathbb{E}\left[ \left\vert \bigtriangleup Y_{t}\right\vert _{H}^{2}\right] +
\mathbb{E}\, \big[ \int_{t}^{T}\left\Vert \bigtriangleup Z_{s}\right\Vert
_{L_{2}\left( E_{1},H\right) }^{2}ds\big] \nonumber \\
& \hspace{1cm}\leq 2\,\mathbb{E}\left[ \int_{t}^{T} \vert  \! \bigtriangleup  \!
\widehat{F}_{s} \vert _{H}\left\vert \bigtriangleup Y_{s}\right\vert _{H}ds\right]
+\mathbb{E}\, \big[ \int_{t}^{T} \Vert \! \bigtriangleup  \! \widehat{G}_{s} 
\Vert_{L_{2}\left( E_{2},H\right)}^{2}\, ds\big] \nonumber \\
&
\hspace{1cm}\leq 2\,C\,\mathbb{E}\left[ \int_{t}^{T}\left( \left\vert
\bigtriangleup Y_{s}\right\vert _{H}+\left\Vert \bigtriangleup
Z_{s}\right\Vert _{L_{2}\left( E_{1},H\right) }\right) \left\vert
\bigtriangleup Y_{s}\right\vert_{H}ds\right] \nonumber \\
&
\hspace{1.7cm}+\, 2 \, C\,\mathbb{E}\left[ \int_{t}^{T}\left( \big{(}\mathbb{E}\, \big[
\left\vert \bigtriangleup Y_{s}\right\vert _{H}^{2}\big] \big{)} ^{\frac{1}{2}}+
\big{(} \mathbb{E}\, \big[ \left\Vert \bigtriangleup Z_{s}\right\Vert _{L_{2}\left(
E_{1},H\right) }^{2}\big] \big{)}^{\frac{1}{2}}\right) \left\vert \bigtriangleup
Y_{s}\right\vert _{H}ds\right] \nonumber \\
&
\hspace{1.7cm}+\,\mathbb{E}\big[ \int_{t}^{T}\left(C \left\vert
\bigtriangleup Y_{s}\right\vert _{H}^{2}+\gamma \, \mathbb{E}\left[ \left\vert
\bigtriangleup Y_{s}\right\vert _{H}^{2}\right] \right) ds\big] \nonumber \\
&
\hspace{1.7cm}+\,\mathbb{E}\, \Big[ \int_{t}^{T}\gamma  \left( \left\Vert
\bigtriangleup Z_{s}\right\Vert _{L_{2}\left( E_{1},H\right) }^{2}+\mathbb{E}
\big[ \left\Vert \bigtriangleup Z_{s}\right\Vert _{L_{2}\left(
E_{1},H\right) }^{2}\big] \right) ds\Big] . &
\end{align}

For any $\varepsilon > 0$, we then observe
\begin{align*}
&  \mathbb{E}\left[ \left \vert \bigtriangleup Y_{t}\right \vert _{H}^{2}\right]
+\mathbb{E}\, \big[ \int_{t}^{T}\left \Vert \bigtriangleup Z_{s}\right \Vert
_{L_2\left( E_{1},H\right) }^{2}ds\big]  \\ &  \hspace{1cm} \leq 2 C \int_{t}^{T}\mathbb{E}\left[
\left( \left \vert \bigtriangleup Y_{s}\right \vert _{H}^{2}+\left \Vert
\bigtriangleup Z_{s}\right \Vert _{L_2\left( E_{1},H\right) }\left \vert
\bigtriangleup Y_{s}\right \vert _{H}\right. \right. \\ &
\hspace{2cm}+  \left. \left. \left \vert \bigtriangleup Y_{s}\right \vert _{H}
\big{(} \mathbb{E}\left[ \left \vert \bigtriangleup Y_{s}\right \vert _{H}^{2}\right] \big{)}^{\frac{1}{2}}
+\left \vert \bigtriangleup Y_{s}\right \vert _{H} \big{(} \mathbb{E}\,
\big[ \left \Vert \bigtriangleup Z_{s}\right \Vert_{L_2\left(
E_{1},H\right)}^{2}\big] \big{)}^{\frac{1}{2}}\right) ds\right] \\ &
\hspace{2cm}+\int_{t}^{T}\left(\left(C+\gamma\right)\, \mathbb{E}\left[ \left \vert \bigtriangleup
Y_{s}\right \vert _{H}^{2}\right] +2 \gamma \, \mathbb{E}\, \big[ \left \Vert
\bigtriangleup Z_{s}\right \Vert_{L_2\left( E_{1},H\right)}^{2}\big]\right) ds
\\ &
\hspace{1cm}\leq 2 C\int_{t}^{T}\left( \left(2+\frac{1}{2\, \varepsilon}+\frac{1}{
2\, \varepsilon}\right) \mathbb{E}\left[ \left \vert \bigtriangleup Y_{s}\right
\vert _{H}^{2}\right] +\left( \frac{\varepsilon}{2}+\frac{\varepsilon}{2}\right)
\mathbb{E}\, \big[ \left \Vert \bigtriangleup Z_{s}\right \Vert _{L_2\left(
E_{1},H\right)}^{2}\big] \right)ds \\ &
\hspace{2cm}+\int_{t}^{T}\left(\left(C+\gamma\right)\, \mathbb{E}\left[ \left \vert \bigtriangleup
Y_{s}\right \vert _{H}^{2}\right] +2\gamma \, \mathbb{E}\, \big[ \left \Vert
\bigtriangleup Z_{s}\right \Vert _{L_2\left( E_{1},H\right)}^{2}\big]\right) ds
\\ &
\hspace{1cm}\leq \int_{t}^{T}\left( \left(5 C+\gamma+\frac{2 C}{\varepsilon}\right) \mathbb{E}
\left[ \left \vert \bigtriangleup Y_{s}\right \vert _{H}^{2}\right] +\left(
2\gamma+2 C\, \varepsilon \right) \mathbb{E}\, \big[ \left \Vert \bigtriangleup
Z_{s}\right \Vert _{L_2\left( E_{1},H\right) }^{2}\big] \right)ds .
\end{align*}
Thus, choosing $\varepsilon=\frac{1-2 \gamma}{4 C}$ (recalling $0< \gamma < 1/2$) yields the following inequality:
\begin{eqnarray*}
&& \hspace{-0.5cm} \mathbb{E}\left[ \left \vert \bigtriangleup Y_{t}\right \vert _{H}^{2}\right]
+\left( \frac{1-2 \gamma}{2}\right) \mathbb{E}\left[ \int_{t}^{T}\left \Vert
\bigtriangleup Z_{s}\right \Vert _{L_2\left( E_{1},H\right) }^{2}ds\right]
\\
&& \hspace{2.2in} \leq\left( 5 C+\gamma+\frac{8C^{2}}{1-2 \gamma}\right) \int_{t}^{T}\mathbb{E}\left[
\left \vert \bigtriangleup Y_{s}\right \vert _{H}^{2}\right] ds.
\end{eqnarray*}
Consequently, by applying Gronwall's inequality, we deduce that $Y_{t}^{1}=Y_{t}^{2}$
and $Z_{t}^{1}=Z_{t}^{2} \;\, a.s.$ for all $0\leq t\leq T.$
\end{proof}

\subsection{Existence of Solutions of MV-FBDSDEs}\label{subs:exist-subsec}
In this section, we establish the existence of solution for the
MV-FBDSDEs~(\ref{eq:2.1}) under assumptions $\left(\mathbf{A}_{1}\right)$--$%
\left(\mathbf{A}_{3}\right)$. We will follow the method of continuation, a method which is explained in \cite{Pe} for the purpose of solving BSDEs with an arbitrary terminal time and also \cite{Yo} for FBSDEs.
\begin{theorem}\label{thm:main thm2}
Under $\left(\mathbf{A}_{1}\right)$--$\left(\mathbf{A}%
_{3}\right)$, MV-FBDSDEs~(\ref{eq:2.1}) has a solution $\left( y,Y,z,Z\right)$ in $\mathfrak{M}^{2}\left( \left[ 0,T\right] ,\mathbb{H}^2 \right) .$
\end{theorem}

We shall employ the method of continuation and divide the proof of this theorem into two separate cases.

\textbf{Case~1:} Let $\theta _{1}>0,\theta _{2}\geq 0$, and $\alpha _{1}>0$. We shall need first Lemma~\ref{lem:lemm2} below which involves a
priori estimates of solutions of the following family of MV-FBDSDEs
parameterized by $\alpha \in \left[ 0,1\right]$:
\begin{equation}\label{eq:3.2}
\left\{
\begin{array}{ll}
dy_{t}=\left( f^{\alpha }\left( t,v_{t},\mathbb{P}_{v_{t}}\right) +\varphi
_{t}\right) dt+\left( g^{\alpha }\left( t,v_{t},\mathbb{P}_{v_{t}}\right) +\phi
_{t}\right)\,  dW_{t}-z_{t}\, \overleftarrow{d B}_{t} , &  \smallskip \\
dY_{t}=\left( F^{\alpha }\left( t,v_{t},\mathbb{P}_{v_{t}}\right) +\psi
_{t}\right) dt+\left( G^{\alpha }\left( t,v_{t},\mathbb{P}_{v_{t}}\right) +\kappa
_{t}\right) \, \overleftarrow{d B}_{t} + Z_{t}\, dW_{t}, &  \smallskip \\
y_{0}=x,\hspace{0.5cm}Y_{T}=h^{\alpha }\left( y_{T},\mathbb{P}
_{y_{T}}\right) +\xi , &
\end{array}%
\right.
\end{equation}%
where $v_t=\left( y_t,Y_t,z_t,Z_t\right) ,$ $\mathbb{P}_{v_t}=\mathbb{P}_{\left(y_t,Y_t,z_t,Z_t\right) },$
$\left( \varphi ,\psi ,\kappa ,\phi \right)
\in \mathfrak{M}^{2}\left( \left[ 0,T\right] ,\mathbb{H}^{2}\right)$
and $\xi \in L^{2}\left( \Omega ,\mathcal{F}_{T},\mathbb{P}%
,H\right) ,$ and for any given $\alpha \in \left[ 0,1\right]:$
\begin{equation*}
\left\{
\begin{array}{ll}
f^{\alpha }\left( t,v_{t},\mathbb{P}_{v_{t}}\right) =\alpha \, f\left(
t,y_{t},Y_{t},z_{t},Z_{t},\mathbb{P}_{v_{t}}\right) , & \smallskip \\
g^{\alpha }\left( t,v_{t},\mathbb{P}_{v_{t}}\right) =\alpha \, g\left(
t,y_{t},Y_{t},z_{t},Z_{t},\mathbb{P}_{v_{t}}\right) , &  \smallskip \\
F^{\alpha }\left( t,v_{t},\mathbb{P}_{v_{t}}\right) =\alpha \,  F\left(
t,y_{t},Y_{t},z_{t},Z_{t},\mathbb{P}_{v_{t}}\right) +\left( 1-\alpha \right) \theta _{1}\left( -y_{t}\right) , &  \smallskip \\
G^{\alpha }\left( t,v_{t},\mathbb{P}_{v_{t}}\right) =\alpha \, G\left(
t,y_{t},Y_{t},z_{t},Z_{t},\mathbb{P}_{v_{t}}\right) +\left( 1-\alpha \right) \theta _{1}\left( -z_{t}\right) , &  \smallskip \\
h^{\alpha }\left( y_{T},P_{y_{T}}\right) =\alpha \, h\left( y_{T},\mathbb{P}%
_{y_{T}}\right) +\left( 1-\alpha \right) y_{T}. &
\end{array}%
\right.
\end{equation*}

When $\alpha =1,$ the existence of the solution of (\ref{eq:3.2}) implies
clearly that of (\ref{eq:2.1}) by letting $\left( \varphi ,\psi ,\phi
,\kappa \right) =\left( 0,0,0,0\right) $. On the other hand, if $\alpha =0$ then (\ref{eq:3.1}) reduces to the following linear FBDSDEs:
\begin{equation}\label{eq:3.3}
\left\{
\begin{array}{ll}
dy_{t}= \varphi _{t}\,  dt+ \phi _{t}\, dW_{t}-z_{t}\,
\overleftarrow{d B}_{t} , &  \smallskip \\
dY_{t}=\left( \theta _{1}\left( -y_{t}\right) +\psi _{t}\right) dt+\left(
\theta _{1}\left( -z_{t}\right) +\kappa _{t}\right)\,  \overleftarrow{d B}_{t}%
+Z_{t}\, dW_{t}, &  \smallskip \\
y_{0}=x,\quad Y_{T}=y_{T}+\xi . &
\end{array}%
\right.
\end{equation}

\begin{lemma}
\label{lem:lemm1} The system~(\ref{eq:3.3}), which is that of (\ref{eq:3.2})
when $\alpha=0,$ has a unique solution $\left( y,Y,z,Z\right) \in \mathfrak{M%
}^{2}\left( \left[ 0,T\right] ,\mathbb{H}^2 \right) .$
\end{lemma}
\begin{proof}
It is easy to verify that MV-FBDSDEs~(\ref{eq:3.3}) satisfies $\left( \mathbf{%
A}_{1}\right) $--$\left( \mathbf{A}_{3}\right) $. From Theorems~(5.3, 5.4)
in \cite{A}, we know that (\ref{eq:3.2}) has a unique solution $\left(
y,Y,z,Z\right) $ in $\mathfrak{M}^{2}\left( \left[ 0,T\right] ,\mathbb{H}%
^{2}\right) $. For more details, we refer the readers to the arguments
presented in \cite{A, AG}.
\end{proof}

\medskip

The following lemma is a key step in the proof of the method of continuation.
\begin{lemma}\label{lem:lemm2}
Assume that $\left(\mathbf{A}_{1}\right)$--$\left(\mathbf{A}_{3}\right)$ holds with
$\theta _{1}>0, \, \theta _{2}\geq 0$, and $\alpha _{1}>0$.
Suppose that there exists a constant $\alpha_{0}\in \left[ 0,1\right)$ such that, for any
$\xi \in L^{2}\left( \Omega ,\mathcal{F}_{T},\mathbb{P},H\right)$ and
 $\left( \varphi ,\psi ,\kappa ,\phi \right)
\in \mathfrak{M}^{2}\left( \left[ 0,T\right] ,\mathbb{H}^{2}\right)$,
MV-FBDSDEs~(\ref{eq:3.2}) has a unique solution.

Then there exists $\delta _{0}\in \left( 0,1\right)$, which only depends  on
$C, \gamma, \alpha_1, \theta_1 , \theta_2$, and $T$,
such that for any $\alpha \in \left[ \alpha _{0},\alpha _{0}+\delta _{0}\right]$,
 MV-FBDSDEs~(\ref{eq:3.2}) has a unique solution.
\end{lemma}
\begin{proof}
Assume that for each $\xi \in L^{2}\left( \Omega ,\mathcal{F}_{T},\mathbb{P}
,H\right)$, $\left( \varphi ,\psi ,\kappa ,\phi \right)
\in \mathfrak{M}^{2}\left( \left[ 0,T\right] ,\mathbb{H}^{2}\right)$,
 MV-FBDSDEs~(\ref{eq:3.2})
has a unique solution for a constant $\alpha =\alpha _{0}\in \left[ 0,1\right)$. Then, for each element $\bar{v}=\left( \bar{y},\bar{Y},\bar{z},\bar{Z}\right)$  of  $\mathfrak{M}^{2}\left(\left[ 0,T\right] ,\mathbb{H}^{2}\right)$, there exists a unique
quadruple $v=\left( y,Y,z,Z\right) \in \mathfrak{M}^{2}\left( \left[ 0,T%
\right] ,\mathbb{H}^{2}\right)$ satisfying the following MV-FBDSDEs:
\begin{eqnarray}\label{eq:new-number}
\left\{
\begin{array}{ll}
dy_{t}=\left( f^{\alpha _{0}}\left( t,v_{t},\mathbb{P}_{v_{t}}\right) +\delta
\, f\left( t,\bar{v}_{t},\mathbb{P}_{\bar{v}_{t}}\right) +\varphi _{t}\right) dt &
\\
\hspace{1.5cm}+\left( g^{\alpha _{0}}\left( t,v_{t},\mathbb{P}_{v_{t}}\right)
+\delta \, g\left( t,\bar{v}_{t},\mathbb{P}_{\bar{v}_{t}}\right) +\phi _{t}\right)
\, dW_{t}-z_{t}\, \overleftarrow{d B}_{t} , &  \medskip \\
dY_{t}=\left( F^{\alpha _{0}}\left( t,v_{t},\mathbb{P}_{v_{t}}\right) +\delta
\left( \theta _{1}\, \bar{y}_{t}
+F\left( t,\bar{v}_{t},\mathbb{P}_{\bar{v}_{t}}\right) \right) +\psi _{t}\right) dt &
 \\
\hspace{1.5cm}+\left( G^{\alpha _{0}}\left( t,v_{t},\mathbb{P}_{v_{t}}\right)
+\delta \left( \theta _{1}\, \bar{z}_{t}
+G\left( t,\bar{v}_{t},\mathbb{P}_{\bar{v}_{t}}\right) \right) +\kappa _{t}\right)\,  \overleftarrow{d B}_{t} + Z_{t}\, dW_{t},
& \medskip \\
y_{0}=x,\hspace{0.5cm}Y_{T}=h^{\alpha _{0}}\left( y_{T},\mathbb{P}%
_{y_{T}}\right) +\delta \left( h\left( \bar{y}_{T},\mathbb{P}_{\bar{y}%
_{T}}\right) -\bar{y}_{T}\right) +\xi . &
\end{array}%
\right.
\end{eqnarray}%

We will show that the mapping defined by
 $I_{\alpha _{0}+\delta }\left( \bar{v},\bar{y}_{T}\right) :=\left( v,y_{T}\right)$
  from  the space $\mathfrak{M}^{2}\left( \left[ 0,T
\right] ,\mathbb{H}^{2}\right) \times L^{2}\left( \Omega ,\mathcal{F}_{T},
\mathbb{P},H\right)$ to itself is a contraction if $\delta >0$ is sufficiently small. To do this, let
$I_{\alpha _{0}+\delta }\left( \bar{v}^{i},\bar{y}_{T}^{i}\right) :=\left(
v^{i},y_{T}^{i}\right) $ for elements $\bar{v}^{i}=\left( \bar{y}^{i},\bar{Y}^{i},\bar{z}^{i},\bar{Z}^{i}\right)$ of $\mathfrak{M}^{2}\left( \left[ 0,T\right] ,\mathbb{H}^{2}\right)$, and let $v^{i,\bar{Z}}=\left( \bar{y}^{i},\bar{Y}^{i},\bar{z}^{i}\right) $ and $v^{i,\bar{z}}=\left( \bar{y}^{i},\bar{Y}
^{i},\bar{Z}^{i}\right) $\ for $i=1,2$. Next, set the notation
$\bigtriangleup \bar{v}=\left( \bigtriangleup \bar{y},\bigtriangleup \bar{Y},\bigtriangleup \bar{z},\bigtriangleup \bar{Z}\right) =\left( \bar{y}^{1}-\bar{y}^{2},\bar{Y}^{1}-\bar{Y}^{2},\bar{z}^{1}-\bar{z}^{2},\bar{Z}^{1}-\bar{Z}^{2}\right)$
and
$\bigtriangleup v=\left( \bigtriangleup y,\bigtriangleup Y,\bigtriangleup
z,\bigtriangleup Z\right) =\left(
y^{1}-y^{2},Y^{1}-Y^{2},z^{1}-z^{2},Z^{1}-Z^{2}\right)$. Also, denote 
 $\bigtriangleup h_{T}=h \big( y_{T}^{1},\mathbb{P}_{y_{T}^{1}} \big)
-h \big( y_{T}^{2},\mathbb{P}_{y_{T}^{2}}\big)$ and $\bigtriangleup \bar{h}_{T}=h
\big( \bar{y}_{T}^{1},\mathbb{P}_{\bar{y}_{T}^{1}}\big) -h \big(
\bar{y}_{T}^{2},\mathbb{P}_{\bar{y}_{T}^{2}} \big)$.

By applying It\^{o}'s formula to $\left\langle \bigtriangleup
y_{t},\bigtriangleup Y_{t}\right\rangle _{H}$, it follows that
\begin{align*}
& \alpha _{0}\, \mathbb{E}\left[\left\langle \bigtriangleup y_{T},\bigtriangleup
h_{T}\right\rangle _{H}\right]+\left( 1-\alpha _{0}\right)\mathbb{E}\left[ \left\vert \bigtriangleup y_{t}\right\vert _{H}^{2} \right]+\delta \,
\mathbb{E}\left[\left\langle \bigtriangleup y_{T},\bigtriangleup \bar{h}_{T}-\bigtriangleup
\bar{y}_{T}\right\rangle _{H}\right] \\ &
\hspace{2.5cm}=\alpha _{0}\,\mathbb{E}\left[ \int_{0}^{T}\left( A\left(
t,v_{t}^{1},\mathbb{P}_{v^1_{t}}\right) -A\left( t,v_{t}^{2},
\mathbb{P}_{v^2_{t}}\right)
,\bigtriangleup v_{t}\right) dt\right] \\ &
\hspace{3cm}-\left( 1-\alpha _{0}\right) \theta _{1}\, \mathbb{E}\left[
\int_{0}^{T}\left(\left\vert \bigtriangleup y_{t}\right\vert _{H}^{2}+\left\Vert
\bigtriangleup z_{t}\right\Vert _{L_{2}\left( E_{2},H\right) }^{2}\right)dt\right]
\\ &
\hspace{3cm}+\,\delta \,\mathbb{E}\left[ \int_{0}^{T}
\left( A\left( t,\bar{v}_{t}^{1},\mathbb{P}_{\bar{v}_{t}^{1}}\right)
 -A\left( t,\bar{v}_{t}^{2},
\mathbb{P}_{\bar{v}_{t}^{2}} \right) ,\bigtriangleup v _{t}\right) dt\right] \\ &
\hspace{3cm}+\,\delta \, \theta _{1}\,\mathbb{E}\left[ \int_{0}^{T}\left(\left\langle
\bigtriangleup y_{t},\bigtriangleup \bar{y}_{t}\right\rangle
_{H}+\left\langle \bigtriangleup z_{t},\bigtriangleup \bar{z}
_{t}\right\rangle _{L_{2}\left( E_{2},H\right) }\right)dt\right] .
\end{align*}
Based on the conditions imposed in $\left(\mathbf{A}_{2}\right)$, we can derive the following inequality:
\begin{align*}
& \left( \alpha _{0}\, \alpha _{1}+\left( 1-\alpha _{0}\right) \right) \mathbb{E}
\left[ \left\vert \bigtriangleup y_{T}\right\vert _{H}^{2}\right] \\ &
\hspace{0.4cm}\leq - \alpha _{0}\,\theta _{2}\, \mathbb{E}\left[ \int_{0}^{T}\left(\left\vert
\bigtriangleup Y_{t}\right\vert _{H}^{2}+\left\Vert \bigtriangleup
Z_{t}\right\Vert _{L_{2}\left( E_{1},H\right) }^{2}\right)dt\right] \\ &
\hspace{0.7cm}-\, \alpha _{0}\,\theta _{1}\, \mathbb{E}\left[ \int_{0}^{T}\left(\left\vert
\bigtriangleup y_{t}\right\vert _{H}^{2}+\left\Vert \bigtriangleup
z_{t}\right\Vert _{L_{2}\left( E_{2},H\right) }^{2}\right)dt\right] \\ &
\hspace{0.7cm}+\, \delta \, \mathbb{E}\left[ \int_{0}^{T}\left\vert A\left( t,\bar{v}
_{t}^{1},\mathbb{P}_{\bar{v}_{t}^{1}}\right) -A\left( t,\bar{v}_{t}^{2},\mathbb{P}_{\bar{v}_{t}^{2}}\right) \right\vert \left\vert \bigtriangleup v_{t}\right\vert dt\right] \\ &
\hspace{0.7cm}+\, \delta \, \theta_{1}\,\mathbb{E}\left[ \int_{0}^{T}
\left( \frac{1}{2}\left\vert \bigtriangleup y_{t}\right\vert_H^{2}+\frac{1}{2}\left\vert
\bigtriangleup \bar{y}_{t}\right\vert_H^{2} + \frac{1}{2}
\left\Vert \bigtriangleup z_{t}\right\Vert_{L_{2}\left( E_{2},H\right)}^{2}+\frac{1}{2}\left\Vert
\bigtriangleup \bar{z}_{t}\right\Vert ^{2}_{L_{2}\left( E_{2},H\right)} \right) dt\right] \\ &
\hspace{0.7cm}+\, \delta \, \mathbb{E}\left[  \frac{1}{2}\left\vert
\bigtriangleup y_{T}\right\vert_H^{2}+\frac{1}{2}\left\vert \bigtriangleup
\bar{y}_{T}\right\vert_H^{2} +\left\vert \bigtriangleup
y_{T}\right\vert_H \left\vert \bigtriangleup \bar{h}_{T}\right\vert_H \right] .
\end{align*}
Therefore, we can rewrite the inequality as follows:
\begin{align*}
&\left( \alpha _{0}\, \alpha _{1}+\left( 1-\alpha _{0}\right)
\right) \mathbb{E}\left[ \left\vert \bigtriangleup y_{T}\right\vert _{H}^{2}
\right] +\alpha _{0}\,\theta _{2}\, \mathbb{E}\left[ \int_{0}^{T}\big{(} \left\vert \bigtriangleup
Y_{t}\right\vert _{H}^{2}+\left\Vert \bigtriangleup Z_{t}\right\Vert
_{L_{2}\left( E_{1},H\right) }^{2}\big{)} \, dt\right] \\ &
\hspace{5.5cm}+\, \alpha _{0}\,\theta _{1}\,\mathbb{E}\left[ \int_{0}^{T}\big{(}\left\vert
\bigtriangleup y_{t}\right\vert _{H}^{2}+\left\Vert \bigtriangleup
z_{t}\right\Vert _{L_{2}\left( E_{2},H\right) }^{2}\big{)} \, dt\right] \\ &
\hspace{0.2cm}\leq \delta \, \mathbb{E}\, \Bigg{[} \int_{0}^{T}\Bigg{(}\frac{1}{2}\left\vert
\bigtriangleup v_{t}\right\vert ^{2}+ C \, \big{(} \left\vert
\bigtriangleup \bar{v}_{t}\right\vert ^{2}+\bar{w}^{2}_{2}\big{(}
\mathbb{P}_{v_{t}^{1}},\mathbb{P}_{v_{t}^{2}}\big{)} \big{)} \\ &
\hspace{3cm}+ \left.\frac{C}{2} \left\vert v_{t}^{1,\bar{Z}}-v_{t}^{2,\bar{Z}}\right\vert^{2}+\frac{\gamma }{2}\, \Big{(}\left\Vert \bigtriangleup \bar{Z}
_{t}\right\Vert _{L_{2}\left( E_{1},H\right) }^{2}+\bar{w}^{2}_{2}\big{(} \mathbb{P}_{\bar{v}_{t}^{1}},\mathbb{P}_{\bar{v}_{t}^{2}}\big{)} \Big{)} \right. \\ &
\hspace{3cm}+  \frac{C}{2} \left\vert v_{t}^{1,\bar{z}}-v_{t}^{2,\bar{z}}\right\vert
^{2}+\frac{\gamma }{2}\, \big{(} \left\Vert \bigtriangleup \bar{z}
_{t}\right\Vert _{L_{2}\left( E_{2},H\right) }^{2}+\bar{w}^{2}_{2}\big{(}\mathbb{P}_{\bar{v}_{t}^{1}},\mathbb{P}_{\bar{v}_{t}^{2}}\big{)} \big{)}\Bigg{)} dt\Bigg{]} \\ &
\hspace{0.2cm} +\, \delta \, \theta _{1}\, \mathbb{E}\, \Bigg{[} \int_{0}^{T}\left( \frac{1}{2}\left( \left\vert \bigtriangleup y_{t}\right\vert _{H}^{2}+\left\vert
\bigtriangleup \bar{y}_{t}\right\vert _{H}^{2}\right) +\frac{1}{2} \big{(}
\left\Vert \bigtriangleup z_{t}\right\Vert _{L_{2}\left( E_{2},H\right)
}^{2}+\left\Vert \bigtriangleup \bar{z}_{t}\right\Vert _{L_{2}\left(
E_{2},H\right) }^{2} \big{)} \right) dt\Bigg{]} \\ &
\hspace{0.2cm} +\, \delta \,\mathbb{E}\, \Big{[} \frac{1}{2}\left\vert
\bigtriangleup y_{T}\right\vert _{H}^{2}+\frac{1}{2}\left\vert
\bigtriangleup \bar{y}_{T}\right\vert _{H}^{2} +\frac{1}{2}\left\vert
\bigtriangleup y_{T}\right\vert _{H}^{2} +\frac{C}{2}\, \Big{(}
\left\vert \bigtriangleup \bar{y}_{T}\right\vert _{H}^{2}+\bar{w}^{2}_{2}\big(
\mathbb{P}_{\bar{y}_{T}^{1}},\mathbb{P}_{\bar{y}_{T}^{2}}\big) \Big{)}
\Big{]} .
\end{align*}

We also know from inequality (\ref{eq:2.2}) that $\bar{w}^{2}_{2}\left(
\mathbb{P}_{v^{1}_t},\mathbb{P}_{v^{2}_t}\right) \leq \mathbb{E}\left[
|\bigtriangleup\!v_t|^{2}\right]$. As a result,
\begin{align*}
& \left( \alpha _{0}\, \alpha_{1}+\left( 1-\alpha _{0}\right) \right) \mathbb{E}
\left[ \left\vert \bigtriangleup y_{T}\right\vert _{H}^{2}\right] +\alpha _{0}\,
\theta_{2}\, \mathbb{E}\left[ \int_{0}^{T} \big{(} \left\vert \bigtriangleup Y_{t}\right\vert
_{H}^{2}+\left\Vert \bigtriangleup Z_{t}\right\Vert _{L_{2}\left(
E_{1},H\right) }^{2} \big{)} \, dt\right] \\ &
\hspace{5.5cm}+ \, \theta _{1}\, \mathbb{E}\left[ \int_{0}^{T} \big{(} \left\vert
\bigtriangleup y_{t}\right\vert _{H}^{2}+\left\Vert \bigtriangleup
z_{t}\right\Vert _{L_{2}\left( E_{2},H\right) }^{2} \big{)} \, dt\right] \\ &
\hspace{2.5cm}\leq \delta\, L\,\mathbb{E}\left[ \int_{0}^{T} \left( \left\vert
\bigtriangleup v _{t}\right\vert ^{2}+\left\vert \bigtriangleup \bar{v}_{t}\right\vert ^{2} \right) \, dt\right] +\delta\, L\, \mathbb{E}\left[ \left\vert
\bigtriangleup y_{T}\right\vert _{H}^{2}+\left\vert \bigtriangleup \bar{y}
_{T}\right\vert _{H}^{2}\right],
\end{align*}
for some generic constant $L>0,$ which from here on  may vary from place to place and depends at most on the constants
$C,$ $\gamma,$ $\alpha_{1} ,$ $\theta _{1},$ $\theta _{2},$ and $T$. Next, since
$$\alpha _{0}\, \alpha_{1}+\left( 1-\alpha _{0}\right) > \min \{1, \alpha_1\} =: \tilde{\alpha}_1$$
 and $\tilde{\alpha}_1 >0$,
then by letting $\beta =\min \left\{ \tilde{\alpha}_1 ,\theta _{1}\right\}$, we deduce that $0< \beta < 1$ and
\begin{align}\label{eq:3.4}
& \mathbb{E}
\left[ \left\vert \bigtriangleup y_{T}\right\vert _{H}^{2}\right] +  \mathbb{E}\left[ \int_{0}^{T} \big{(} \left\vert
\bigtriangleup y_{t}\right\vert _{H}^{2}+\left\Vert \bigtriangleup
z_{t}\right\Vert _{L_{2}\left( E_{2},H\right) }^{2} \big{)} \, dt\right]  \nonumber \\
& \hspace{1cm}\leq \frac{\delta \, L}{\beta} \left(\mathbb{E}\left[ \int_{0}^{T} \big{(} \left\vert
\bigtriangleup v _{t}\right\vert ^{2}+\left\vert \bigtriangleup \bar{v}
_{t}\right\vert ^{2} \big{)} \, dt\right] + \mathbb{E}\left[ \left\vert
\bigtriangleup y_{T}\right\vert _{H}^{2}+\left\vert \bigtriangleup \bar{y}
_{T}\right\vert _{H}^{2}\right] \right),
\end{align}
after neglecting the term
$\alpha _{0}\,\theta_{2}\, \mathbb{E}\left[ \int_{0}^{T} \big{(} \left\vert \bigtriangleup Y_{t}\right\vert
_{H}^{2}+\left\Vert \bigtriangleup Z_{t}\right\Vert _{L_{2}\left(
E_{1},H\right) }^{2} \big{)} \, dt\right]$ that contains $\alpha_0$, as we aim to find a value of $\delta$ that is independent of $\alpha_0$.

We therefore need to find estimates for
$\mathbb{E}\left[ \int_{0}^{T} \big{(} \left\vert
\bigtriangleup Y_{t}\right\vert _{H}^{2}+\left\Vert \bigtriangleup
Z_{t}\right\Vert _{L_{2}\left( E_{1},H\right) }^{2} \big{)} \, dt\right]$. To this end, we apply It\^{o}'s formula to $\left\vert \bigtriangleup
Y_{t}\right\vert _{H}^{2}$ and take the expectation. Eventually, we find that
\begin{align*}
& \mathbb{E}\left[ \left\vert \bigtriangleup Y_{t}\right\vert _{H}^{2}\right] + \mathbb{E} \Big{[} \int_{t}^{T}\left\Vert \bigtriangleup Z_{s}\right\Vert
_{L_{2}\left( E_{1},H\right) }^{2}ds\Big{]}=\mathbb{E}\left[ \left\vert \bigtriangleup
h^{\alpha _{0}}\left( y_{T},\mathbb{P}_{y_{T}}\right) +\delta \left(
\bigtriangleup \bar{h}_{T}-\bigtriangleup \bar{y}_{T}\right) \right\vert _{H}^{2}\right]
\nonumber \\ &
\hspace{2cm} - \, 2\, \mathbb{E} \Big{[} \int_{t}^{T}\Big{(} \left\langle \,  \bigtriangleup
F^{\alpha _{0}}\left( s,v_{s},\mathbb{P}_{v_{s}}\right) ,\bigtriangleup
Y_{s}\right\rangle _{H} \nonumber \\ &
\hspace{6cm}
 + \left\langle \, \delta \left( \theta _{1}\left(
\bigtriangleup \bar{y}_{s}\right) +\bigtriangleup F\left( s,\bar{v}_{s},\mathbb{P}_{\bar{v}_{s}}\right) \right) ,\bigtriangleup Y_{s}\right\rangle _{H}\Big{)} \, ds%
\Big{]}
\nonumber \\ &
\hspace{1.5cm} +\,  \mathbb{E}\, \Big{[} \int_{t}^{T}\left\Vert \bigtriangleup
G^{\alpha _{0}}\left( s,v_{s},\mathbb{P}_{v_{s}}\right)  +\delta \left( \theta _{1}\left( \bigtriangleup
\bar{z}_{s}\right) +\bigtriangleup G\left( s,\bar{v}_{s},
\mathbb{P}_{\bar{v}_{s}}\right) \right) \right\Vert _{L_{2}\left( E_{2},H\right) }^{2}\, ds\Big{]},
\end{align*}
where%
\begin{equation*}
\left\{
\begin{array}{ll}
\bigtriangleup F^{\alpha _{0}}\left( s,v_{s},\mathbb{P}_{v_{s}}\right) =\alpha
_{0}\bigtriangleup\!F\left( s,v_{s},\mathbb{P}_{\left(
y_{s},Y_{s},z_{s},Z_{s}\right) }\right) +\left( 1-\alpha _{0}\right) \theta
_{1}\left( -\left( \bigtriangleup y_{s}\right) \right) , &  \smallskip \\
\bigtriangleup G^{\alpha _{0}}\left( s,v_{s},\mathbb{P}_{v_{s}}\right) =\alpha
_{0}\bigtriangleup\!G\left( s,v_{s},\mathbb{P}_{\left(
y_{s},Y_{s},z_{s},Z_{s}\right) }\right) +\left( 1-\alpha _{0}\right) \theta
_{1}\left( -\left( \bigtriangleup z_{s}\right) \right) , &  \smallskip \\
\bigtriangleup h^{\alpha _{0}}\left( y_{T},\mathbb{P}_{y_{T}}\right) =\alpha
_{0}\bigtriangleup\!h_{T} +\left( 1-\alpha
_{0}\right) \left( \bigtriangleup y_{T}\right) , &  \smallskip \\
\bigtriangleup F\left( s,\bar{v}_{s},\mathbb{P}_{\bar{v}_{s}}\right)
 =F\left( t,\bar{v}_{s}^{1},\mathbb{P}_{\left( \bar{y}_{s}^{1},\bar{Y}_{s}^{1},\bar{z}
_{s}^{1},\bar{Z}_{s}^{1}\right) }\right) -F\left( t,\bar{v}_{s}^{2},
\mathbb{P}_{\left( \bar{y}_{s}^{2},\bar{Y}_{s}^{2},\bar{z}_{s}^{2},\bar{Z}
_{s}^{2}\right) }\right) , &  \smallskip \\
\bigtriangleup G\left( s,\bar{v}_{s},\mathbb{P}_{\bar{v}_{s}}\right)
=G\left( t,\bar{v}_{s}^{1},\mathbb{P}_{\left( \bar{y}_{s}^{1},
\bar{Y}_{s}^{1},\bar{z}_{s}^{1},\bar{Z}_{s}^{1}\right) }\right) -G\left( t,\bar{v}_{s}^{2},\mathbb{P}_{\left( \bar{y}_{s}^{2},
\bar{Y}_{s}^{2},\bar{z}_{s}^{2},\bar{Z}_{s}^{2}\right) }\right) . &
\end{array}
\right.
\end{equation*}
It follows that
\begin{equation}\label{eq:3.5}
\mathbb{E}\left[ \left\vert \bigtriangleup Y_{t}\right\vert _{H}^{2}\right] +
\mathbb{E}\left[ \int_{t}^{T}\left\Vert \bigtriangleup Z_{s}\right\Vert
_{L_{2}\left( E_{1},H\right) }^{2}ds\right]  \leq I_{1}+I_{2} (t) +I_{3} (t) +I_{4}(t),
\end{equation}
where
\begin{align*}
& I_1=   4\, \mathbb{E}\left[ \alpha _{0}^{2}\left\vert
\bigtriangleup h_{T} \right\vert
_{H}^{2}+\left( 1-\alpha _{0}\right) ^{2}\left\vert \bigtriangleup
y_{T} \right\vert _{H}^{2}
+  \delta ^{2}\left\vert \bigtriangleup \bar{h}_{T} \right\vert _{H}^{2}+\delta ^{2}\left\vert
 \bigtriangleup \bar{y}_{T} \right\vert _{H}^{2}\right] , \smallskip \\
& I_2(t) :=  2\, \mathbb{E}\left[ \int_{t}^{T}\Big{(} \left\vert \bigtriangleup
F^{\alpha _{0}}\left( s,v_{s},\mathbb{P}_{v_{s}}\right) \right\vert _{H}\left\vert
\bigtriangleup Y_{s}\right\vert _{H}\right.  \\ & \hspace{3.5cm} \left.
 +\, \left( \delta \, \theta _{1}\left\vert
\bigtriangleup \bar{y}_{s}\right\vert _{H}+\left\vert \delta
\bigtriangleup\!F\left( s,\bar{v}_{s},\mathbb{P}_{\bar{v}_{s}}\right) \right\vert
_{H}\right) \left\vert \bigtriangleup Y_{s}\right\vert _{H} \Big{)} \, ds\right] ,
\smallskip \\
& I_3(t) := \mathbb{E}\left[ \int_{t}^{T} \left(\frac{1+\gamma }{2\gamma }\right) \alpha
_{0}^{2}\left\Vert \bigtriangleup G\left( s,v_{s},\mathbb{P}_{\left(
y_{s},Y_{s},z_{s},Z_{s}\right) }\right) \right\Vert _{L_{2}\left(
E_{2},H\right) }^{2} \, ds \right] ,
\smallskip \\
& I_4(t) := 3\, \left(\frac{1+\gamma }{1-\gamma }\right) \mathbb{E}\Bigg[ \int_{t}^{T} \Big{(}  \left( 1-\alpha _{0}\right)
^{2}\theta _{1}^{2}\left\Vert \bigtriangleup z_{s}\right\Vert _{L_{2}\left(
E_{2},H\right) }^{2} + \delta ^{2}\, \theta
_{1}^{2}\left\Vert \bigtriangleup \bar{z}_{s}\right\Vert _{L_{2}\left(
E_{2},H\right) }^{2} \\ & \hspace{3in} 
+ \, \delta ^{2}\left\Vert \bigtriangleup G\left( s,\bar{v}
_{s},\mathbb{P}_{\bar{v}_{s}}\right) \right\Vert _{L_{2}\left( E_{2},H\right)
}^{2} \Big{)} \, ds\Bigg] .
\end{align*}

On the other hand, from the fact that $h$ is a Lipschitz mapping, we obtain
\begin{equation}\label{eq:3.6}
I_1 \leq C\, \mathbb{E}\left[ \left\vert \bigtriangleup
y_{T} \right\vert _{H}^{2}
+  \delta \left\vert
 \bigtriangleup \bar{y}_{T} \right\vert _{H}^{2}\right] .
\end{equation}
Since $F$ is Lipschitz, we also have
\begin{align*}
\hspace{-0.4cm} I_2(t) \leq \, \mathbb{E}\left[ \int_{t}^{T}\Bigg{(} \Big{(}\frac{16\, C\, \alpha_{0}^{2}}{1-\gamma }\Big{)}
\left\vert \bigtriangleup Y_{s}\right\vert _{H}^{2}+ \Big{(}\frac{1-\gamma}{16\, C}\Big{)} \, C \, \Big( \left\vert \bigtriangleup v_{s}\right\vert ^{2}
+ \bar{w}^{2}_{2}\left(\mathbb{P}_{v^{1}_t},\mathbb{P}_{v^{2}_t}\right) \Big) \right.&  \\
\hspace{1cm} \left.
+\, (1-\alpha_{0})\, \theta_{1}\left(\left\vert \bigtriangleup Y_{s}\right\vert _{H}^{2}+\left\vert \bigtriangleup y_{s}\right\vert _{H}^{2} \right)+\, \delta \, \theta _{1}\left( \left\vert
\bigtriangleup \bar{y}_{s}\right\vert _{H}^{2}+\left\vert \bigtriangleup Y_{s}\right\vert _{H}^{2}\right)  \right. & \\
\hspace{2cm} + \, \delta \, \bigg( \left\vert \bigtriangleup Y_{s}\right\vert _{H}^{2}+C\, \Big( \left\vert \bigtriangleup \bar{v}_{s}\right\vert ^{2}
+ \bar{w}^{2}_{2}\left(\mathbb{P}_{\bar{v}^{1}_t},\mathbb{P}_{\bar{v}^{2}_t}\right) \Big) \bigg{)} \Bigg)\, ds & \Bigg]. &
\end{align*}
Thus, by using inequality (\ref{eq:2.2}), it follows that
\begin{align}\label{eq:3.7}
\hspace{-0.2cm} I_2(t) \leq C\, \mathbb{E}\left[ \int_{t}^{T}\Big{(}
\left\vert \bigtriangleup Y_{s}\right\vert _{H}^{2}+ \left\vert \bigtriangleup y_{s}\right\vert _{H}^{2}
+ \left\Vert \bigtriangleup z_{s}\right \Vert_{L_2\left(  E_{2},H\right)}^{2} \Big{)} \, ds\right] + C \, \delta\, \mathbb{E}\left[ \int_{t}^{T}\left\vert \bigtriangleup \bar{y}_{s}\right\vert _{H}^{2}ds\right]& \nonumber \\
\hspace{-2.5cm} + \, \Big{(}\frac{1-\gamma}{8}\Big{)} \,\mathbb{E}\left[ \int_{t}^{T}\left\Vert \bigtriangleup Z_{s}\right \Vert_{L_2\left(  E_{1},H\right)}^{2} ds\right]+C\, \delta\, \mathbb{E}\left[ \int_{t}^{T}\left\vert \bigtriangleup \bar{v}_{s}\right\vert ^{2}ds\right]. &
\end{align}

For $I_3(t)$ and $I_4(t)$, we apply $\left(\mathbf{A}_{1}\right)$ to see that
\begin{align}\label{eq:3.8}
 & I_3(t) \nonumber \\ & \leq \, \mathbb{E}\left[ \int_{t}^{T} \Big{(} \frac{1+\gamma}{4\gamma}\Big{)}\, \alpha_{0}\,\Big{(}
C\left\vert v^{1,Z}-v^{2,Z}\right\vert ^{2}+ \gamma\left(\left\Vert \bigtriangleup Z_{s}\right \Vert_{L_2\left(  E_{1},H\right)}^{2}
+ \bar{w}^{2}_{2}\left(\mathbb{P}_{v^{1}_t},\mathbb{P}_{v^{2}_t}\right)\right)\Big{)} \, ds\right]  \nonumber \\
& \leq \,C\, \mathbb{E}\left[ \int_{t}^{T}\Big{(}
\left\vert \bigtriangleup Y_{s}\right\vert _{H}^{2}+ \left\vert \bigtriangleup y_{s}\right\vert _{H}^{2}
+ \left\Vert \bigtriangleup z_{s}\right \Vert_{L_2\left(  E_{2},H\right)}^{2} \Big{)} \, ds\right] 
\nonumber \\
&
\hspace{2in} + \, \Big{(}\frac{1+\gamma}{2}\Big{)}\, \alpha_{0}\,\mathbb{E}\left[ \int_{t}^{T}\left\Vert \bigtriangleup Z_{s}\right \Vert_{L_2\left(  E_{1},H\right)}^{2} ds\right]
\end{align}
and
\begin{align}\label{eq:3.9-a}
 I_4(t) & \leq 3 \, \Big{(}\frac{4\gamma}{1-\gamma}\Big{)}\, \mathbb{E}\, \Bigg{[} \int_{t}^{T}\Bigg{(}(1-\alpha_{0})^{2}\, \theta_{1}^{2}\,
\left\Vert \bigtriangleup z_{s} \right \Vert_{L_2\left(  E_{2},H\right)}^{2} + \delta^{2}\, \theta_{1}^{2} \, \left\Vert \bigtriangleup \bar{z}_{s}
\right \Vert_{L_2\left(  E_{2},H\right)}^{2}
\nonumber \\ & \hspace{1.5cm}
 +\, \delta^{2} \bigg( C \left\Vert v^{1,\bar{z}}-v^{2,\bar{z}}\right \Vert_{L_2\left(  E_{2},H\right)}^{2} + \gamma \, \Big(\left\vert \bigtriangleup \bar{z}_{s}\right\vert _{H}^{2}
+ \bar{w}^{2}_{2}\left(\mathbb{P}_{\bar{v}^{1}_t},\mathbb{P}_{\bar{v}^{2}_t}\right)\Big)\bigg)\Bigg{)} \, ds\Bigg{]} \nonumber \\ &
 \leq \,C\, \mathbb{E}\Bigg[ \int_{t}^{T}\Big{(}
\left\Vert \bigtriangleup z_{s}\right \Vert_{L_2\left(  E_{2},H\right)}^{2} +  \delta \left\vert \bigtriangleup \bar{y}_{s}\right\vert _{H}^{2}
+ \delta\left\vert \bigtriangleup \bar{Y}_{s}\right\vert _{H}^{2}+ \delta\left\Vert \bigtriangleup \bar{z}_{s}\right \Vert_{L_2\left(  E_{2},H\right)}^{2}
 \nonumber \\
& \hspace{2in}
 +\, \delta\, C\, \Big{(}\frac{8\gamma^{2}}{1-\gamma}\Big{)} \left\Vert \bigtriangleup \bar{Z}_{s}\right \Vert_{L_2\left(  E_{1},H\right)}^{2} \Big{)} \, ds \Bigg].
\end{align}

Now, substitute (\ref{eq:3.6})--(\ref{eq:3.9-a}) into (\ref{eq:3.5}) and use $\alpha_{0}<1$ to find that there exists a universal constant $L>0$, which is, of course, independent of $\alpha_0$, such that
\begin{eqnarray}\label{eq:3.10}
&& \hspace{-1cm}\mathbb{E}\left[ \left\vert \bigtriangleup Y_{t}\right\vert _{H}^{2}\right] +\left(1- \frac{5+3\gamma}{8}\right)
\mathbb{E}\left[ \int_{t}^{T}\left\Vert \bigtriangleup Z_{s}\right\Vert
_{L_{2}\left( E_{1},H\right) }^{2}ds\right] \nonumber \\ &&
 \leq \,C\, \mathbb{E}\left[ \int_{t}^{T} \left\vert \bigtriangleup Y_{s}\right\vert
_{H}^{2} \, ds\right]+L\, \mathbb{E}\left[ \left( \left\vert
\bigtriangleup y_{T}\right\vert _{H}^{2}+\delta \left\vert \bigtriangleup
\bar{y}_{T}\right\vert _{H}^{2}\right) \right] \nonumber \\  &&
\hspace{2cm}+\,L\,\mathbb{E}\left[ \int_{t}^{T}\left( \left\vert
\bigtriangleup y_{s}\right\vert _{H}^{2}+\left\Vert \bigtriangleup
z_{s}\right\Vert _{L_{2}\left( E_{2},H\right) }^{2}+\delta \left\vert
\bigtriangleup \bar{v}_{s}\right\vert ^{2}\right) ds\right].
\end{eqnarray}
Recall that $0<\gamma<1/2$ and apply Gronwall's inequality to obtain
\begin{eqnarray*}
&& \hspace{-1cm}\mathbb{E}\left[ \left\vert \bigtriangleup Y_{t}\right\vert _{H}^{2}\right]
 \leq \,L \, e^{C\, (T-t)}\, \Bigg{(} L\, \mathbb{E}\left[ \left\vert
\bigtriangleup y_{T}\right\vert _{H}^{2}+\delta \left\vert \bigtriangleup
\bar{y}_{T}\right\vert _{H}^{2} \right] \nonumber \\  &&
\hspace{2.5cm} +\,L\,\mathbb{E}\, \left[ \int_{t}^{T}\left( \left\vert
\bigtriangleup y_{s}\right\vert _{H}^{2}+\left\Vert \bigtriangleup
z_{s}\right\Vert _{L_{2}\left( E_{2},H\right) }^{2}+\delta \left\vert
\bigtriangleup \bar{v}_{s}\right\vert ^{2}\right) ds\right]\Bigg{)},
\end{eqnarray*}
for each $t\in[0,T]$, although $0< \gamma <1$ is only needed here.  As a result, we deduce
\begin{align}\label{eq:3.11}
& \mathbb{E}\left[ \int_{0}^{T}\big{(} \left\vert \bigtriangleup Y_{s}\right\vert
_{H}^{2}+\left\Vert \bigtriangleup Z_{s}\right\Vert _{L_{2}\left(
E_{1},H\right) }^{2} \big{)}  \, ds\right] \leq L\, \mathbb{E}\left[ \left\vert
\bigtriangleup y_{T}\right\vert _{H}^{2}+\delta \left\vert \bigtriangleup
\bar{y}_{T}\right\vert _{H}^{2} \right]   \nonumber \\ &
\hspace{2.4cm}+\,L\,\mathbb{E}\left[ \int_{0}^{T}\left( \left\vert
\bigtriangleup y_{s}\right\vert _{H}^{2}+\left\Vert \bigtriangleup
z_{s}\right\Vert _{L_{2}\left( E_{2},H\right) }^{2}+\delta \left\vert
\bigtriangleup \bar{v}_{s}\right\vert ^{2}\right) ds\right] ,
\end{align}
recalling that the constant $L$ can vary from one place to another.

Next, we combine  the crucial results (\ref{eq:3.4}) and (\ref{eq:3.11}) to obtain
\begin{align*}
& \hspace{-0.3cm}  \mathbb{E}\left[ \left\vert \bigtriangleup y_{T}\right\vert
_{H}^{2}\right] + \mathbb{E}\left[ \int_{0}^{T}\left\vert \bigtriangleup v
_{t}\right\vert ^{2}dt\right] \\
& \hspace{1cm} \leq \frac{\delta\, L}{\beta} \left( \mathbb{E}\left[ \int_{0}^{T} \big{(} \left\vert
\bigtriangleup v _{t}\right\vert ^{2}+\left\vert \bigtriangleup \bar{v}%
_{t}\right\vert ^{2} \big{)} \, dt\right] + \mathbb{E}\left[ \left\vert
\bigtriangleup y_{T}\right\vert _{H}^{2}+\left\vert \bigtriangleup \bar{y}%
_{T}\right\vert _{H}^{2}\right] \right) \\ & \hspace{1.7cm}
+\, L\, \mathbb{E}\left[ \left\vert
\bigtriangleup y_{T}\right\vert _{H}^{2} \right] + L\,\mathbb{E}\left[ \int_{0}^{T}\left( \left\vert
\bigtriangleup y_{t}\right\vert _{H}^{2}+\left\Vert \bigtriangleup
z_{t}\right\Vert _{L_{2}\left( E_{2},H\right) }^{2} \right) \, dt \right] .
\end{align*}%
Then we apply (\ref{eq:3.4}) once more to the last two terms of this latter inequality to derive
\begin{align*}
& \hspace{-0.3cm} \mathbb{E}\left[ \int_{0}^{T}\left\vert \bigtriangleup v
_{t}\right\vert ^{2}dt\right]  + \mathbb{E}\left[ \left\vert \bigtriangleup y_{T}\right\vert
_{H}^{2}\right] \\
& \hspace{1cm} \leq \frac{\delta\, L}{\beta} \left( \mathbb{E}\left[ \int_{0}^{T} \big{(} \left\vert
\bigtriangleup v _{t}\right\vert ^{2}+\left\vert \bigtriangleup \bar{v}%
_{t}\right\vert ^{2} \big{)} \, dt\right] + \mathbb{E}\left[ \left\vert
\bigtriangleup y_{T}\right\vert _{H}^{2}+\left\vert \bigtriangleup \bar{y}%
_{T}\right\vert _{H}^{2}\right] \right) .
\end{align*}%
By taking $\delta\leq \delta_0:=\frac{\beta}{3L}$, it follows that
\begin{equation*}
\begin{array}{ll}
\mathbb{E}\left[ \int_{0}^{T}\left \vert \bigtriangleup v_{t}\right \vert
^{2}dt+\left \vert \bigtriangleup y_{T}\right \vert _{H}^{2}\right] \leq
\frac{1}{2}\, \mathbb{E}\left[ \int_{0}^{T}\left \vert \bigtriangleup \bar{v}%
_{t}\right \vert ^{2}dt+\left \vert \bigtriangleup \bar{y}_{T}\right \vert
_{H}^{2}\right] . &
\end{array}%
\end{equation*}
Hence, we conclude that the mapping $I_{\alpha _{0}+\delta
} $ is a contraction for all fixed $\delta$ in $[0,\delta_0]$. As a result, 
$I_{\alpha _{0}+\delta}$  attains a unique fixed point $\left(y,Y,z,Z\right)$ in
 $\mathfrak{M}^{2}\left( \left[ 0,T\right] ,\mathbb{H}^{2}\right) ,$ which is the solution of MV-FBDSDE~(\ref{eq:3.2}) for $\alpha =\alpha
_{0}+\delta , \, \delta \in [0,\delta_0]$.
\end{proof}

\bigskip

\textbf{Case~2:} Let $\theta _{1}\geq 0, \, \theta _{2}>0$, and $\alpha _{1}\geq 0$. Consider the following family of MV-FBDSDEs,
parameterized by $\alpha \in \left[ 0,1\right]$:
\begin{eqnarray}\label{eq:3.14}
\left \{
\begin{array}{ll}
dy_{t}=\left( \breve{f}^{\alpha }\left( t,v_{t},\mathbb{P}_{v_{t}}\right) +\varphi
_{t}\right) dt+\left( \breve{g}^{\alpha }\left( t,v_{t},\mathbb{P}_{v_{t}}\right) +\phi
_{t}\right) \, dW_{t}-z_{t}\, \overleftarrow{d B}_{t} , \smallskip  \\
dY_{t}=\left( \breve{F}^{\alpha }\left( t,v_{t},\mathbb{P}_{v_{t}}\right) +\psi
_{t}\right) dt+\left( \breve{G}^{\alpha }\left( t,v_{t},\mathbb{P}_{v_{t}}\right) +\kappa
_{t}\right)  \overleftarrow{d B}_{t}+Z_{t}\, dW_{t},   \smallskip  \\
 y_{0}=x,\hspace{0.5cm}Y_{T}=\breve{h}^{\alpha }\left( y_{T},\mathbb{P}
_{y_{T}}\right) +\xi ,
\end{array}
 \right.
\end{eqnarray}%
where
\begin{equation*}
\left\{
\begin{array}{ll}
\breve{f}^{\alpha }\left( t,v_{t},\mathbb{P}_{v_{t}}\right) =\alpha\, f\left(
t,y_{t},Y_{t},z_{t},Z_{t},\mathbb{P}_{\left( y_{t},Y_{t},z_{t},Z_{t}\right)
}\right) +\left( 1-\alpha \right) \theta _{2}\left( -Y_{t}\right) ,&  \smallskip \\
\breve{g}^{\alpha }\left( t,v_{t},\mathbb{P}_{v_{t}}\right) =\alpha \, g\left(
t,y_{t},Y_{t},z_{t},Z_{t},\mathbb{P}_{\left( y_{t},Y_{t},z_{t},Z_{t}\right)
}\right) +\left( 1-\alpha \right) \theta _{2}\left( -Z_{t}\right) , &  \smallskip \\
F^{\alpha }\left( t,v_{t},\mathbb{P}_{v_{t}}\right) =\alpha \, F\left(
t,y_{t},Y_{t},z_{t},Z_{t},\mathbb{P}_{\left( y_{t},Y_{t},z_{t},Z_{t}\right)
}\right) , &  \smallskip \\
\breve{G}^{\alpha }\left( t,v_{t},\mathbb{P}_{v_{t}}\right) =\alpha \, G\left(
t,y_{t},Y_{t},z_{t},Z_{t},\mathbb{P}_{\left( y_{t},Y_{t},z_{t},Z_{t}\right)
}\right) , &  \smallskip \\
\breve{h}^{\alpha }\left( y_{T},P_{y_{T}}\right) =\alpha \, h\left( y_{T},\mathbb{P}
_{y_{T}}\right) . &
\end{array}
\right.
\end{equation*}

When $\alpha =1,$ the existence of the solution of (\ref{eq:3.14}) clearly implies
 that of (\ref{eq:2.1}) by letting
 $\left( \varphi ,\psi ,\phi ,\kappa \right) =\left( 0,0,0,0\right)$. On the other hand, if $\alpha =0,$ (\ref{eq:3.14}) is uniquely solvable as explained in Case~1. We now state a crucial lemma that will help us complete the proof of Theorem~\ref{thm:main thm2}.
\begin{lemma}
\label{lem:lemm3}
Assume that $\left(\mathbf{A}_{1}\right)$--$\left(
\mathbf{A}_{3}\right)$ hold with $\theta _{1}\geq 0,\, \theta _{2}>0$, and $\alpha
_{1}\geq 0$. Suppose there exists a constant $\alpha
_{0}\in \left[ 0,1\right) $ such that, for any
$\xi \in L^{2}\left( \Omega ,\mathcal{F}_{T},\mathbb{P},H\right)$ and
$\left( \varphi ,\psi ,\kappa ,\phi \right)
\in \mathfrak{M}^{2}\left( \left[ 0,T\right] ,\mathbb{H}^{2}\right)$,
MV-FBDSDEs~(\ref{eq:3.14}) has a unique solution.

Then there exists $\delta _{0}\in \left( 0,1\right) $ which only depends  on
$C, \gamma, \alpha_1, \theta_1, \theta _{2}$, and $T$,  such that for any $\alpha
\in \left[ \alpha _{0},\alpha _{0}+\delta _{0}\right]$, MV-FBDSDEs~(\ref{eq:3.14}) has a unique solution.
\end{lemma}
\begin{proof}
Assume that, for all $\xi \in L^{2}\left( \Omega ,\mathcal{F}_{T},\mathbb{P}
,H\right)$ and $\left( \varphi ,\psi ,\kappa ,\phi \right)
\in \mathfrak{M}^{2}\left( \left[ 0,T\right] ,\mathbb{H}^{2}\right)$,
MV-FBDSDEs~(\ref{eq:3.14})
has a unique solution for a constant $\alpha =\alpha _{0}\in \left[ 0,1\right)$. Then, for each element $\bar{v}=\left( \bar{y},\bar{Y},\bar{z},\bar{Z}\right)$ of
 $\mathfrak{M}^{2}\left(\left[ 0,T\right] ,\mathbb{H}^{2}\right)$,
 there exists a unique quadruple $v=\left( y,Y,z,Z\right) \in \mathfrak{M}^{2}\left( \left[ 0,T\right] ,\mathbb{H}^{2}\right)$ satisfying the following MV-FBDSDEs:
\begin{equation*}
\left\{
\begin{array}{ll}
dy_{t}=\left( \breve{f}^{\alpha _{0}}\left( t,v_{t},\mathbb{P}_{v_{t}}\right) +\delta\,
\left( \theta _{2}\, \bar{Y}_{t}+f\left( t,\bar{v}_{t},\mathbb{P}_{\bar{v}_{t}}\right) \right)
 +\varphi _{t}\right) dt &  \\
\hspace{2.5cm}+\left( g^{\alpha _{0}}\left( t,v_{t},\mathbb{P}_{v_{t}}\right)
+\delta \left( \theta _{2}\, \bar{Z}_{t}
+g\left( t,\bar{v}_{t},\mathbb{P}_{\bar{v}_{t}}\right) \right) +\phi _{t}\right)  dW_{t}-z_{t}\, \overleftarrow{d B}_{t} , &
\medskip \\
dY_{t}=\left( \breve{F}^{\alpha _{0}}\left( t,v_{t},\mathbb{P}_{v_{t}}\right) +\delta\,
F\left( t,\bar{v}_{t},\mathbb{P}_{\bar{v}_{t}}\right) +\psi _{t}\right) dt &  \\
\hspace{2.5cm}+\left( \breve{G}^{\alpha _{0}}\left( t,v_{t},\mathbb{P}_{v_{t}}\right)
+\delta \, G\left( t,\bar{v}_{t},\mathbb{P}_{\bar{v}_{t}}\right) +\kappa _{t}\right)
\overleftarrow{d B}_{t}+Z_{t}\, dW_{t}, & \medskip \\
y_{0}=x,\hspace{0.5cm}Y_{T}=\breve{h}^{\alpha _{0}}\left( y_{T},\mathbb{P}
_{y_{T}}\right) +\delta \,  h\left( \bar{y}_{T},\mathbb{P}_{\bar{y}_{T}}\right)
+\xi . &
\end{array}
\right.
\end{equation*}

We argue as in Case~1.
Let us consider the mapping $I_{\alpha_{0}+\delta}$ defined in the proof of
Lemma~\ref{lem:lemm1} and retain the same notations as set there after system~(\ref{eq:new-number}). By applying It\^{o}'s formula to
$\left\langle \bigtriangleup y_{t},\bigtriangleup Y_{t}\right\rangle_{H}$ and disregarding the terms involving $\alpha_0$, we obtain
\begin{align}\label{eq:3.15}
&  \theta_{2}\, \mathbb{E}\left[ \int_{0}^{T}  \left\vert
\bigtriangleup v_{t}\right\vert _{H}^{2} dt\right] \leq \delta \, L\,\mathbb{E}\left[ \int_{0}^{T} \left( \left\vert
\bigtriangleup v _{t}\right\vert ^{2}+ \left\vert \bigtriangleup \bar{v}%
_{t}\right\vert ^{2}\right) dt\right] \nonumber \\ &
\hspace{2.2in}  +\, \delta \, \mathbb{E}\left[ \left\vert
\bigtriangleup y_{T}\right\vert _{H}^{2} \right] + \delta \, L\, \mathbb{E}\left[ \left\vert \bigtriangleup \bar{y}_{T}\right\vert _{H}^{2}\right] .
\end{align}

On the other hand, we can follow a similar approach as in
(\ref{eq:3.5}) by employing It\^{o}'s formula for
$\left\vert \bigtriangleup y_{s}\right\vert_{H}^{2}$ to get
\begin{equation}\label{eq:3.16}
\mathbb{E}\left[ \left\vert
\bigtriangleup y_{T}\right\vert _{H}^{2}\right] + \mathbb{E}\left[ \int_{0}^{T}  \left\Vert
\bigtriangleup z_{t}\right\Vert_{L_2(E_2,H)}^{2} dt\right] \leq  L\,\mathbb{E}\left[ \int_{0}^{T} \left( \left\vert
\bigtriangleup v _{t}\right\vert ^{2}+ \delta \left\vert \bigtriangleup \bar{v}
_{t}\right\vert ^{2}\right) dt\right] .
\end{equation}
These two inequalities play a crucial role here.

Now, let $\beta':=\min\{\theta_2, 1\}$ to observe that $0 < \beta' < 1$. By combining (\ref{eq:3.15}) and (\ref{eq:3.16}), we can derive the following inequality:
\begin{align*}
& \beta' \left( \mathbb{E}\left[ \left\vert \bigtriangleup y_{T}\right\vert
_{H}^{2}\right] +\mathbb{E}\left[ \int_{0}^{T}\left\vert \bigtriangleup v
_{t}\right\vert ^{2}dt\right] \right) \\ &
\hspace{1cm} \leq  \delta \, L\,\mathbb{E}\left[ \int_{0}^{T} \left( \left\vert
\bigtriangleup v _{t}\right\vert ^{2}+ \left\vert \bigtriangleup \bar{v}
_{t}\right\vert ^{2}\right) dt\right] + \delta \, L\, \mathbb{E}\left[ \left\vert
\bigtriangleup y_{T}\right\vert _{H}^{2} \right] + \delta \, L\, \mathbb{E}\left[ \left\vert \bigtriangleup \bar{y}_{T}\right\vert _{H}^{2}\right]
\\ &
\hspace{2.8in} + \, L\,\mathbb{E}\left[ \int_{0}^{T} \left( \left\vert
\bigtriangleup v _{t}\right\vert ^{2}+ \delta \left\vert \bigtriangleup \bar{v}
_{t}\right\vert ^{2}\right) dt\right] .
\end{align*}
Furthermore, applying (\ref{eq:3.15}) again to the term $\mathbb{E}\left[ \int_{0}^{T} \left\vert
\bigtriangleup v _{t}\right\vert ^{2} dt\right]$ and utilizing the preceding inequality, we obtain
\begin{align*}
& \mathbb{E}\left[ \int_{0}^{T}\left\vert \bigtriangleup v
_{t}\right\vert ^{2}dt + \left\vert \bigtriangleup y_{T}\right\vert
_{H}^{2}\right]   \\
& \hspace{1cm} \leq  \frac{\delta \, L}{\beta'} \left( \mathbb{E}\left[ \int_{0}^{T} \left( \left\vert
\bigtriangleup v _{t}\right\vert ^{2}+ \left\vert \bigtriangleup \bar{v}%
_{t}\right\vert ^{2}\right) dt\right] +  \mathbb{E}\left[ \left\vert
\bigtriangleup y_{T}\right\vert _{H}^{2} \right] +  \mathbb{E}\left[ \left\vert \bigtriangleup \bar{y}_{T}\right\vert _{H}^{2}\right]\right) .
\end{align*}%
Thus, if we choose $\delta\leq\delta_0 :=\frac{\beta'}{3 L}$,  we conclude that
\begin{equation*}
\begin{array}{ll}
\mathbb{E}\left[ \int_{0}^{T}\left \vert \bigtriangleup v_{t}\right \vert
^{2}dt+\left \vert \bigtriangleup y_{T}\right \vert _{H}^{2}\right] \leq
\frac{1}{2}\, \mathbb{E}\left[ \int_{0}^{T}\left \vert \bigtriangleup \bar{v}%
_{t}\right \vert ^{2}dt+\left \vert \bigtriangleup \bar{y}_{T}\right \vert
_{H}^{2}\right] . &
\end{array}
\end{equation*}

The remainder of the proof follows a similar approach as in Lemma~\ref{lem:lemm2}.
\end{proof}

\bigskip

We emphasize that the condition $0< \gamma <1/2$ in assumption $\left(\mathbf{A}_{1}\right)$ is needed frankly in Case~2 to establish the proof of
Theorem~\ref{thm:main thm1}
and so the proof of the preceding lemma.
The reader can find similar details in \cite[Lemma~3.8]{AG2}.

\bigskip

We are now ready to conclude the proof of Theorem~\ref{thm:main thm2}.

\medskip

\noindent\begin{proof}[Proof completion of Theorem~\protect\ref{thm:main thm2}]
In Case~1 (when $\theta _{2}>0$), we already know that for each element $\xi $ of
$L^{2}\left( \Omega ,\mathcal{F}_{T},\mathbb{P};H\right) $ and
$\left(\varphi ,\psi ,\phi ,\kappa \right) \in \mathfrak{M}^{2}\left( \left[ 0,T\right],\mathbb{H}^{2}\right)$, the MV-FBDSDEs~(\ref{eq:3.2}) has a unique solution when $\alpha =0$. It then follows from Lemma~\ref{lem:lemm2} that there exists a positive constant
 $\delta _{0}=\delta_{0}\left( C,\gamma ,\alpha_1 ,\theta _{1}, \theta _{2}, T\right)$ such that for any $\delta \in \left[ 0,\delta _{0}\right] ,$ $\xi \in L^{2}\left( \Omega ,\mathcal{F}_{T},\mathbb{P};H\right)$, and $\left( \varphi ,\psi ,\phi
,\kappa \right) \in \mathfrak{M}^{2}\left( \left[ 0,T\right] ,\mathbb{H}^{2}\right)$, (\ref{eq:3.2}) has a unique solution for $\alpha =\delta $. Moreover, since $\delta _{0}$\ depends only on $C,\gamma ,\alpha_1 ,\theta _{1}, \theta_2, T$, we can repeat this process $N$\ times with ${1\leq N\delta_{0}<1+\delta _{0}}$. In particular, for $\alpha =1$\ with
 $\left( \varphi,\psi ,\phi ,\kappa \right) \equiv 0,\phi \equiv 0$, we deduce that
MV-FBDSDEs~(\ref{eq:2.1}) has a unique solution in
 $\mathfrak{M}^{2}\left( \left[ 0,T\right] ,\mathbb{H}^{2}\right)$.

\medskip

For Case~2 (when $\alpha _{1}>0$ and $\theta _{1}>0$), given any
$\xi \in L^{2}\left( \Omega ,\mathcal{F}_{T},\mathbb{P};H\right)$ and $\left( \varphi
,\psi ,\phi ,\kappa \right) \in \mathfrak{M}^{2}\left( \left[ 0,T\right] ,\mathbb{H}^{2}\right)$, MV-FBDSDEs~(\ref{eq:3.14}) has a unique solution when $\alpha =0$. Consequently, Lemma~\ref{lem:lemm2} implies that there exists a constant $\delta_{0} >0$ that depends only on $C,\gamma , \alpha_1, \theta_1, \theta _{2}$, and $T$, such that, for any element
 $\delta \in \left[ 0,\delta _{0}\right]$, $\xi \in L^{2}\left( \Omega ,\mathcal{F}_{T},\mathbb{P};H\right)$, and $\left( \varphi ,\psi ,\phi ,\kappa \right) \in
\mathfrak{M}^{2}\left( \left[ 0,T\right] ,\mathbb{H}^{2}\right)$, system~(\ref{eq:3.2}) has a unique solution for $\alpha =\delta $. Therefore, similar to the preceding case, we conclude that the MV-FBDSDEs~(\ref{eq:2.1}) attains a unique solution in
 $\mathfrak{M}^{2}\left( \left[ 0,T\right] ,\mathbb{H}^{2}\right)$.
\end{proof}

\medskip

We conclude this section by providing two examples to illustrate the results of
Theorems~(\ref{thm:main thm1}, \ref{thm:main thm2}) and to demonstrate how to handle our conditions.
\begin{example}\label{ex:ex1}
Let $E$ and $H$ be two real separable Hilbert spaces. Suppose $B$ and $W$ are cylindrical Wiener processes on $E$. Consider the following system on $H$:
\begin{equation}\label{eq:3.9}
\left\{
\begin{array}{ll}
dy_{t}=\left( \frac{1}{2}\, \mathbb{E}\left[ Y_{t}\right] -Y_{t}\right) dt+\left( \frac{1}{4}\,\mathbb{E}\left[ Z_{t}\right] - \frac{1}{2} \, Z_{t}\right) dW_{t}-z_{t}\, \overleftarrow{d B}_{t} , & \smallskip \\
dY_{t}=\left(\frac{1}{2}\,\mathbb{E}\left[ y_{t}\right] -y_{t}\right) dt
+ \left( \frac{1}{4}\,\mathbb{E}\left[ z_{t}\right] - \frac{1}{2}\, z_{t}\right) \overleftarrow{d B}_{t}+Z_{t}\, dW_{t}, & \smallskip \\
y_{0}=x \; (\in H),\quad Y_{T}=-\frac{1}{2}\,\mathbb{E}\left[ y_{T}\right] + y_{T} .
&
\end{array}%
\right.
\end{equation}
In order to relate this system to MV-FBDSDEs~(\ref{eq:2.1}), we define for $t\in [0,T]$,
\begin{equation*}
\begin{array}{ll}
f\left( t,y_{t},Y_{t},z_{t},Z_{t},\mathbb{P}_{\left( y_{t},Y_{t},z_{t},Z_{t}\right)
}\right)=\frac{1}{2}\, \mathbb{E}\left[ Y_{t}\right] -Y_{t}, &  \smallskip \\
g\left( t,y_{t},Y_{t},z_{t},Z_{t},\mathbb{P}_{\left( y_{t},Y_{t},z_{t},Z_{t}\right)
}\right)=\frac{1}{4}\,\mathbb{E}\left[ Z_{t}\right] -\frac{1}{2} \, Z_{t}, & \smallskip
\\
F\left( t,y_{t},Y_{t},z_{t},Z_{t},\mathbb{P}_{\left( y_{t},Y_{t},z_{t},Z_{t}\right)
}\right)=\frac{1}{2}\,\mathbb{E}\left[ y_{t}\right] -y_{t}, &  \smallskip \\
G\left( t,y_{t},Y_{t},z_{t},Z_{t},\mathbb{P}_{\left( y_{t},Y_{t},z_{t},Z_{t}\right)
}\right)=\frac{1}{4}\,\mathbb{E}\left[ z_{t}\right] - \frac{1}{2}\, z_{t}, & \smallskip
 \\
h\left( y_{T},\mathbb{P}_{y_{T}}\right) =-\frac{1}{2}\,\mathbb{E}\left[ y_{T}\right] + y_{T}.
\end{array}
\end{equation*}
In particular, we have used here
$$\mathbb{E}\left[ Y_{t}\right]= \int_H x_1 \, d \mathbb{P}_{Y_t} (x_1) = \int_{\mathbb{H}^2} \Psi (x_1,x_2,x_3,x_4) \, d \mathbb{P}_{\left( y_{t},Y_{t},z_{t},Z_{t}\right)},$$
through Fubini's theorem, where
$\Psi (x_1,x_2,x_3,x_4)=x_1\cdot 1 \cdot 1 \cdot 1 = x_1$. Similar expressions hold for  $\mathbb{E}\left[ Z_{t}\right]$, $\mathbb{E}\left[ y_{t}\right]$, and $\mathbb{E}\left[ z_{t}\right]$. So the dependence of these mappings $f,g,F,G,h$ on a measure $\mu\in \mathcal{P}_2(\mathbb{H}^2)$ is only through its first moment $\int u \, d \mu$.

\medskip

Now, with the help of (\ref{eq:2.2}) and the Cauchy-Schwarz inequality, we observe
\begin{equation*}
\begin{array}{ll}
\left\vert f\left( t,v_t^{1},\mathbb{P}_{v_t^{1}}\right) -f\left( t,v_t^{2},
\mathbb{P}_{v_t^{2}}\right) \right\vert \leq \left\vert Y_{t}^{1}-Y_{t}^{2}\right\vert +\frac{1}{2}\, \bar{w}
_{2}\left( \mathbb{P}_{Y_{t}^{1}},\mathbb{P}_{Y_{t}^{2}}\right) ,& \smallskip \\
\left\vert F\left( t,v_t^{1},\mathbb{P}_{v_t^{1}}\right) -F\left( t,v_t^{2},
\mathbb{P}_{v_t^{2}}\right) \right\vert \leq \left\vert
y_{t}^{1}-y_{t}^{2}\right\vert + \frac{1}{2}\, \bar{w}_{2}\left( \mathbb{P}_{y_{t}^{1}},
\mathbb{P}_{y_{t}^{2}}\right) ,& \smallskip \\
\left\vert g\left( t,v_t^{1},\mathbb{P}_{v_t^{1}}\right) -g\left( t,v_t^{2},
\mathbb{P}_{v_t^{2}}\right) \right\vert ^{2}\leq \frac{1}{2} \left\Vert Z_{t}^{1}-Z_{t}^{2}\right\Vert
^{2}+\frac{1}{8}\, \bar{w}^{2}_{2}\left( \mathbb{P}_{Z_{t}^{1}},\mathbb{P}_{Z_{t}^{2}}\right),
 & \smallskip \\
\left\vert G\left( t,v_t^{1},\mathbb{P}_{v_t^{1}}\right) -G\left( t,v_t^{2},
\mathbb{P}_{v_t^{2}}\right) \right\vert ^{2}\leq  \frac{1}{2} \left\Vert z_{t}^{1}-z_{t}^{2}\right\Vert
^{2}+\frac{1}{8}\, \bar{w}^{2}_{2}\left( \mathbb{P}_{z_{t}^{1}},\mathbb{P}_{z_{t}^{2}}\right). &
\end{array}
\end{equation*}
We also have
\begin{align*}
& \mathbb{E} \left[ \left( A\left( t,v_t^{1},\mathbb{P}_{v_t^{1}}\right) -A\left(
t,v_t^{2},\mathbb{P}_{v_t^{2}}\right) ,v_t^{1}-v_t^{2}\right)  \right] \\ &
\hspace{2.5cm}
\leq - \frac{1}{4} \, \mathbb{E} \left[ \, \left\vert \Delta
y_{t}\right\vert^{2} + \left\vert
\Delta Y_{t}\right\vert ^{2}+\left\Vert\Delta z_{t}\right\Vert ^{2} +\left\Vert \Delta Z_{t}\right\Vert^{2} \, \right]
\end{align*}
and
\begin{equation*}
\mathbb{E} \left[ \left\langle h\left( y_{T}^{1},\mathbb{P}%
_{y_{T}^{1}}\right) -h\left( y_{T}^{2},\mathbb{P}_{y_{T}^{2}}\right)
,y_{T}^{1}-y_{T}^{2}\right\rangle \right] \geq
 \frac{1}{2}\, \mathbb{E} \left[ \left\vert \Delta y_{T}\right\vert^{2} \, \right] .
\end{equation*}
Therefore, by setting $ C=1, \gamma=\frac{1}{8} \, , \,  \theta_{1}=\theta _{2}=\frac{1}{4} \,,$ and $\alpha_1=\frac{1}{2} \,,$ it follows that assumptions
 $\left(\mathbf{A}_{1}\right) $--$\left(\mathbf{A}_{3}\right)$ are satisfied. As a result, based on Theorems~(\ref{thm:main thm1}, \ref{thm:main thm2}), we deduce that
system~(\ref{eq:3.9}) has a unique solution.
\end{example}

We will now provide a counter example to show that the assumption
$\left(\mathbf{A}_{2}\right)$ in Theorems~(\ref{thm:main thm1}, \ref{thm:main thm2}) is necessary and cannot be dropped.
\begin{example}\label{ex:ex2}
Let us consider the following MV-FBDSDEs on $H=\mathbb{R}$, with spaces
$E_1=E_2=\mathbb{R}:$
\begin{equation}\label{eq:3.10-ex}
\left\{
\begin{array}{ll}
dy_{t}=\mathbb{E}\left[ Y_{t}\right] dt -z_{t}\,
\overleftarrow{d B}_{t} , &  \smallskip \\
dY_{t}=- \, \mathbb{E}\left[ y_{t}\right]  dt-z_{t}\, \overleftarrow{d B}_{t}+Z_{t}\, dW_{t}, & \smallskip \\
y_{0}=0,\hspace{0.5cm}Y_{T}=-\, \mathbb{E}\left[ y_{T}\right] , &
\end{array}
\right.
\end{equation}
for $ T= \frac{3 \pi}{4}$. Here, $B$ and $W$ are $1$-dimensional Brownian motions.

Using the notation set in assumption $\left(\mathbf{A}_{2}\right)$, we have, for
 $v=\left( y,Y,z,Z\right)$,
\begin{equation*}
\begin{array}{ll}
A\left( t,v_t,\mathbb{P}_{v_t}\right) =\left( - \, \mathbb{E}\left[
y_{t}\right]  , \mathbb{E}\left[
Y_{t}\right] ,-z_{t},0\right)  . &
\end{array}
\end{equation*}
Moreover, noting that
\begin{align*}
& \mathbb{E} \left[ \left( A\left( t,v_t^{1},\mathbb{P}_{v_t^{1}}\right) -A\left(
t,v_t^{2},\mathbb{P}_{v_t^{2}}\right) ,v_t^{1}-v_t^{2}\right) \right] \\
& \hspace{1cm} = \, - \left(
 \mathbb{E}\left[\Delta y_{t}\right] \right) ^{2} + \left(  \mathbb{E}\left[\Delta Y_{t}\right] \right) ^{2}-\mathbb{E}\left[ \left\Vert \left( \Delta z_{t}\right)\right \Vert^{2}\right],
\end{align*}
we realize that assumption $\left(\mathbf{A}_{2}\right)$ does not hold.
As a result, (\ref{eq:3.10-ex}) might not have a unique solution. Indeed,
$(\sin t, \cos t, 0,0)$ is  a solution of (\ref{eq:3.10-ex}) in addition to the trivial solution $(y_t,Y_t,z_t,Z_t)=(0,0,0,0)$.
\end{example}

\section{Application to Stochastic Optimal Control}\label{sec:4}
Let $E$ be a separable Hilbert space, and let $U$ be a nonempty convex subset of $K$. We say that $u_{\cdot}:[0,T]\times \Omega \rightarrow K $
is \emph{admissible} if $u_{\cdot}\in \mathfrak{M}^{2}\left( \left[ 0,T\right]
,K \right)$ and $u_{t}\in U$ for each $t\in [0,T]$. The set of all such admissible controls will be denoted by $\mathcal{U}_{ad}$. In this section, we establish sufficient optimality conditions for a stochastic control problem governed by MV-FBDSDEs over infinite-dimensional separable real Hilbert spaces.
In particular, the stochastic control problem we consider aims to minimize the  \emph{cost functional} (or \emph{objective functional}):
\begin{eqnarray}\label{eq:4.1}
 && J\left( u_{\cdot}\right) =\mathbb{E}\big[\varphi \left( y_{T}^{u_{\cdot}},
\mathbb{P}_{y_{T}^{u_{\cdot}}}\right) +\psi \left( Y_{0}^{u_{\cdot}},
\mathbb{P}_{Y_{0}^{u_{\cdot}}}\right) \nonumber \\
&& \hspace{2.5cm} +\int_{0}^{T}\ell \big(
t,y_{t}^{u_{\cdot}},Y_{t}^{u_{\cdot}},z_{t}^{u_{\cdot}},Z_{t}^{u_{\cdot
}},u_{t },\mathbb{P}_{\left( y_{t}^{u_{\cdot}},Y_{t}^{u_{\cdot
}},z_{t}^{u_{\cdot}},Z_{t}^{u_{\cdot}}\right) }\big)dt\big]  
\end{eqnarray}
over $\mathcal{U}_{ad}$, subject to the state dynamic:
\begin{equation}\label{eq:4.2}
\left \{
\begin{array}[c]{l}
dy_{t}^{u_{\cdot}}=f\left( t,y_{t}^{u_{\cdot}},Y_{t}^{u_{\cdot
}},z_{t}^{u_{\cdot}},Z_{t}^{u_{\cdot}},u_{t },\mathbb{P}_{\left(
y_{t}^{u_{\cdot}},Y_{t}^{u_{\cdot}},z_{t}^{u_{\cdot}},Z_{t}^{u_{\cdot
}}\right) }\right) dt   \\
\hspace{2.5cm} + \,g\left( t,y_{t}^{u_{\cdot}},Y_{t}^{u_{\cdot
}},z_{t}^{u_{\cdot}},Z_{t}^{u_{\cdot}},u_{t },\mathbb{P}_{\left(
y_{t}^{u_{\cdot}},Y_{t}^{u_{\cdot}},z_{t}^{u_{\cdot}},Z_{t}^{u_{\cdot
}}\right) }\right) dW_{t}-z_{t}^{u_{\cdot}}\, d\overleftarrow{B_{t}} , \medskip \\
dY_{t}^{u_{\cdot}}=-F\left( t,y_{t}^{u_{\cdot}},Y_{t}^{u_{\cdot
}},z_{t}^{u_{\cdot}},Z_{t}^{u_{\cdot}},u_{t },\mathbb{P}_{\left(
y_{t}^{u_{\cdot}},Y_{t}^{u_{\cdot}},z_{t}^{u_{\cdot}},Z_{t}^{u_{\cdot
}}\right) }\right) dt   \\
\hspace{2.5cm} - \, G\left( t,y_{t}^{u_{\cdot}},Y_{t}^{u_{\cdot
}},z_{t}^{u_{\cdot}},Z_{t}^{u_{\cdot}},u_{t },\mathbb{P}_{\left(
y_{t}^{u_{\cdot}},Y_{t}^{u_{\cdot}},z_{t}^{u_{\cdot}},Z_{t}^{u_{\cdot
}}\right) }\right) d\overleftarrow{B_{t}}+Z_{t}^{u_{\cdot}}\, dW_{t}, \medskip  \\
y_{0}^{u_{\cdot}}=x,\text{ }Y_{T}^{u_{\cdot}}= c\, y_{T}^{u_{\cdot}} + \xi ,
\end{array}
\right.
\end{equation}
with coefficients: 
\begin{center}
  $\left(  f,F\right)  :[0,T]\times \mathbb{H}^2 \times K \times \mathcal{P}_{2}\left( \mathbb{H}^2 \right)  \rightarrow H$,  \\ $g:[0,T]\times\mathbb{H}^2 \times K \times \mathcal{P}_{2}\left( \mathbb{H}^2 \right)
\rightarrow L_2\left(  E_{1},H\right)$,   \\
$G:[0,T]\times\mathbb{H}^2 \times K \times \mathcal{P}_{2}\left( \mathbb{H}^2 \right)
\rightarrow L_2\left(  E_{2},H\right)$,
  \end{center}  
being measurable mappings so that (\ref{eq:4.1}) is defined,  $\xi$ is an $\mathcal{F}_T$-measurable random variable, and $c$ is a constant.

We say that $u^{\ast}_{\cdot}\in \mathcal{U}_{ad}$
is an \emph{optimal control} if it satisfies
\begin{equation}\label{eq:4.3}
J(u^{\ast}_{\cdot })=\inf_{u_{\cdot}\in \mathcal{U}_{ad}}J(u_{\cdot}).
\end{equation}

To address this control problem (\ref{eq:4.1})--(\ref{eq:4.3}), we need to introduce the concept of \emph{$L$-differentiability} with respect to
probability measure. This is necessary due to the
dependence of distribution appearing in both (\ref{eq:4.1}) and (\ref{eq:4.2}). We can then obtain the adjoint equations of (\ref{eq:4.2}), which resemble the MVDSDEs studied in Section~\ref{sec:3}.

In our control problem (\ref{eq:4.1})--(\ref{eq:4.3}), both the state process and the cost functional depend on the distribution $\mathbb{P}_{\left(y_{t}^{u_{\cdot}},Y_{t}^{u_{\cdot}},z_{t}^{u_{\cdot}},Z_{t}^{u_{\cdot}}\right)}$ of the state process, providing more generality to cover cases such as those considered in Examples~(\ref{ex:ex1}, \ref{ex:ex2}).

\subsection{The $L$-Differentiability and Convexity of Functions of Measures}\label{sec:appl}
In this subsection, we recall the definition of the $L$-derivative of functions of measures. The $L$-derivative was introduced by P. Lions, and in this regard, we refer to \cite[Chapter~5]{CD2} for more details on such a notion. Bensoussan et al., \cite{Bens023}, gave an alternative equivalent definition. We shall be working over Hilbert spaces. The idea is to view the probability measures in $\mathcal{P}_{2}\left(E\right)$ over a separable real Hilbert space $E$ as laws of random variables 
$X \in L^{2}\left(\Omega ,\mathcal{F},\mathbb{P},E\right)$ so that $\mu =\mathbb{P}_{X}$. To be more precise, we assume that probability space 
$\left(\Omega ,\mathcal{F},\mathbb{P}\right)$ is rich enough in the sense that for every $\mu \in \mathcal{P}_{2}\left(E\right)$, there is a random variable $X \in L^{2}\left( \Omega ,\mathcal{F},\mathbb{P},E\right)$ such that $\mu =\mathbb{P}_{X}$. A function $\Phi :\mathcal{P}_{2}\left(E\right) \to \mathbb{R}$ is said to be 
\emph{$L$-differentiable} at $\mu_0$ if there exists 
$X\in L^{2}\left( \Omega ,\mathcal{F},
\mathbb{P},E\right)$ with $\mu_0 =\mathbb{P}_{X_0}$ such that the \emph{lifted function} 
$\hat{\Phi}:L^{2}\left( \Omega ,\mathcal{F},
\mathbb{P},E\right) \to \mathbb{R}$, given by $\hat{\Phi}\left(X\right) :=\Phi \left(\mathbb{P}_{X}\right)$ for $X \in L^{2}\left( \Omega ,\mathcal{F},\mathbb{P},E\right)$, is Fr\'{e}chet differentiable at $X_0$, i.e. there exists a continuous linear functional
$D\hat{\Phi} \left(X_0\right) : L^{2}\left( \Omega ,\mathcal{F},\mathbb{P},E\right) \rightarrow \mathbb{R}$ satisfying
$$\hat{\Phi}(X_0 +\triangle  X) - \hat{\Phi}(X_0)=D\hat{\Phi}(X_0)(\triangle  X)
+o(\Vert \triangle  X \Vert),$$
where $\triangle X$ represents a perturbation.

By the Riesz representation theorem, there exists a unique random variable 
$\zeta_0$ in $L^{2}\left( \Omega ,\mathcal{F},\mathbb{P},E\right)$ such that
$D\hat{\Phi} \left(X_0\right) \left(X\right) = \langle
 \zeta_0 , X \rangle$, for each
 $X \in L^{2}\left( \Omega ,\mathcal{F},\mathbb{P}
,E\right)$, where $\langle \cdot, \cdot \rangle$ denotes the inner product in 
$L^{2}\left( \Omega ,\mathcal{F},\mathbb{P},E\right)$. 

It is known (see  \cite{Cardaliaguet} and \cite{CD2}) that there exists a measurable function $\rho : H \to H$ depending only on $\mu_0$ such that $\zeta_0 = \rho(Y )\; \, a.s.$ for all $Y$ with $\mathbb{P}_Y = \mu_0$. We define the \emph{$L$-derivative} 
$\partial_{\mu}\Phi \left(\mu_0 \right)(Y)$ of 
$\hat{\Phi}$ at $\mu_0$ along the
random variable $Y$ by $\rho(Y)$. Therefore, we have $a.s.$
 $\partial_{\mu}\Phi \left(\mathbb{P}_Y \right)(Y)= \rho(Y)= \nabla \hat{\Phi}\left(X_0\right)$,
where $\nabla \hat{\Phi}\left(X_0\right)$ is the gradient of $\hat{\Phi}$ at the point $X_0$.

The continuity of $\partial _{\mu }\Phi \left( x,\mu
\right) $ is understood as the continuity of the mapping 
$X \mapsto \partial _{\mu }\Phi \left(\mathbb{P}_{X}\right)
 \left( X\right)$ from  $L^{2}\left( \Omega ,\mathcal{F},\mathbb{P}
,E\right)$ to $L^{2}\left( \Omega ,\mathcal{F},\mathbb{P},E\right)$.

Similarly, for each fixed $t\in \left[ 0,T\right]$, a
function $\Phi :\mathcal{P}_{2}\left( \mathbb{H}^2\right) \to \mathbb{R}$ is \emph{differentiable} at $\mu$ if there exists a
quadruple of random variables $\left(y,Y,z,Z\right)$ in $\mathbb{H}^2$ with 
$\mu =\mathbb{P}_{\left( y,Y,z,Z\right) }$ so that the lifted function $\hat{\Phi}$, given
by $\hat{\Phi}\left(y,Y,z,Z\right) =\Phi \left(\mathbb{P}_{\left(y,Y,z,Z\right) }\right)$, is Fr\'{e}chet differentiable at $\left(y,Y,z,Z\right)$. The partial 
$L$-derivatives $\partial _{\mu
_{y}}\Phi ,$ $\partial _{\mu_{_Y}}\Phi ,$ $\partial _{\mu_{_z}}\Phi$, and $\partial _{\mu_{_Z}}\Phi $ at ${\mu}$ along $\left( y,Y,z,Z\right)$ 
can be viewed uniquely as an element $\nabla \hat{\Phi}\left(V\right)$ of
$L^{2}\left( \Omega ,\mathcal{F},\mathbb{P},\mathbb{H}^2\right)$, which can be represented as
$\left( \partial _{\mu_{_y}}\Phi \left(\mathbb{P}_{V}\right)  ,\partial _{\mu_{_Y}}\Phi \left(\mathbb{P}_{V}\right)  ,\partial _{\mu_{_z}}\Phi \left(\mathbb{P}_{V}\right),\partial _{\mu_{_Z}}\Phi \left(\mathbb{P}_{V}\right) \right) \left( V\right) ,$
where $V=\left( y,Y,z,Z\right)$.

Finally, let us introduce the following notation. Consider $( \tilde{\Omega},
\mathcal{\tilde{F}},\tilde{\mathbb{P}})$ as a copy of the probability space
 $( \Omega ,\mathcal{F},\mathbb{P})$. For any pair of random variables $(X, X')$ in $L^{2}(\Omega, \mathcal{F}, \mathbb{P}, E) \times L^{2}(\Omega, \mathcal{F}, \mathbb{P}, E)$, we denote their independent copies on $(\tilde{\Omega}, \mathcal{\tilde{F}}, \tilde{\mathbb{P}})$ as $(\tilde{X}, \tilde{X'})$. Furthermore, we denote the expectation under the probability measure $\tilde{\mathbb{P}}$ as $\mathbb{\tilde{E}}$. 

We say that $\Phi$ is \emph{$L$-convex} (or merely \emph{convex}) if for every $\mu , \mu' \in \mathcal{P}_2(E)$, we have
\begin{equation}\label{eq:convexity}
\Phi(\mu') -  \Phi(\mu) - \tilde{\mathbb{E}} \left[\left \langle \, \partial_{\mu} \Phi (\mu) (X) , X'-X  \right \rangle \right] \geq 0 ,
\end{equation}
whenever $X, X' \in L^2(\tilde{\Omega}, \mathcal{\tilde{F}}, \tilde{\mathbb{P}})$ with distributions $\mu$ and $\mu'$, respectively.

\subsection{The Maxmum Principle}\label{subsec:4.2}
To establish the maximum principle for optimality, we need the following assumptions.

$\left(\mathbf{A}_{4}\right)$: Assume that
\begin{equation*}
\left\{
\begin{array}{ll}
\hspace{-0.15cm} {\rm (i)} & \hspace{-0.25cm} F,f,G,g,\ell \;
\text{are continuous and continuously Fr\'{e}chet differentiable with }    \\
& \hspace{-0.25cm} \text{respect to} \; (y,Y,z,Z,u) \in \mathbb{H}^2 \times K, \text{ and}\;\varphi ,\psi \;
\text{are continuously differentiable} \\ 
& \hspace{-0.25cm}  \text{with respect to} \; y\in H \; \text{and}\;Y\in H,\, \text{respectively}. \medskip \\
\hspace{-0.15cm} {\rm (ii)} & \hspace{-0.25cm} \text{The Fr\'{e}chet derivatives of }
F,f,G,g\; \text{with respect to the above arguments} \\
& \hspace{-0.25cm} \text{are continuous and bounded, uniformly in} \; \left( t,\mu
\right). \; \text{Moreover, the Fr\'{e}chet}  \\
& \hspace{-0.25cm} \text{derivatives of}\;  \phi =g,G \; \text{satisfy } \left\vert \frac{\partial}{\partial z} \phi \left(
t,y,Y,z,Z,v,\mu \right) \right\vert^2 <\gamma  \text{ and}  \\
& \hspace{-0.25cm} \left\vert \frac{\partial}{\partial Z} \phi \left( t,y,Y,z,Z,v,\mu \right)
\right\vert^2 <\gamma , \text{ with } 0<\gamma<\frac{1}{6}. 
\medskip \\
\hspace{-0.2cm} {\rm (iii)} & \hspace{-0.25cm} \text{The derivatives of} \; \ell
\text{\ are bounded by } \\
& \text{ }C\, (1+\left\vert y\right\vert _{H}+\left\vert Y\right\vert
_{H}+\left\Vert z\right\Vert _{L_{2}\left( E_{2},H\right) }+\left\Vert
Z\right\Vert _{L_{2}\left( E_{1},H\right) }+\bar{w}_{2}\left( \mu
,\delta _{0}\right) ). \medskip \\
\hspace{-0.15cm} {\rm (iv)} & \hspace{-0.25cm} \text{The derivatives of}\;\varphi
\;\text{and}\;\psi \;\text{are bounded by}\; C\, (1+\left\vert y\right\vert_{H}+\bar{w}_{2}\left(\mu ,\delta _{0}\right) )\;\text{and} \\
& \hspace{-0.2cm} C\, (1+\left\vert Y\right\vert_{H}+\bar{w}_{2}\left( \mu ,\delta _{0}\right) ), \;\text{respectively},
\end{array}
\right.
\end{equation*}
for some constant $C>0$, where 
$\delta_{0}$ denotes the Dirac measure at $0$.

\bigskip

$\left( \mathbf{A}_{5}\right)$: Suppose that the following conditions hold:
\begin{equation*}
\left\{
\begin{array}{ll}
\hspace{-0.15cm} {\rm (i)} & \hspace{-0.25cm} F,f,G,g,\ell \;
\text{are } L\text{-differentiable with respect to } \mu \in \mathcal{P}_{2}\left( \mathbb{H}^{2}\right), \text{and}\;\varphi ,\psi \;\text{are} \\
& \hspace{-0.25cm} \text{continuously } L\text{-differentiable
with} \text{ respect to}\;\mu . \medskip \\
\hspace{-0.15cm} {\rm (ii)} & \hspace{-0.25cm}\text{The}\; 
L\text{-derivatives of }F,f,G,g \text{ are continuous and bounded, uniformly in } \\
& \hspace{-0.25cm} \left(t,y,Y,z,Z,v,\mu \right);
 \text{in particular, we require} \; \medskip \\
& \displaystyle \int_{\mathbb{H}^{2}}\Big\vert \partial_{\mu_{_z}} \phi \left( t,y,Y,z,Z,v,\mu \right) \big(y',Y',z',Z'\big) \Big\vert ^{2}d\mu \big(y',Y',z',Z'\big) <\frac{1}{3} \, \gamma , \; \text{and} \medskip \\
& \displaystyle \int_{\mathbb{H}^{2}}\Big\vert \partial_{\mu_{_Z}} \phi \left(
t,y,Y,z,Z,v,\mu \right)\big(y',Y',z',Z'\big) \Big\vert ^{2}d\mu \big(y',Y',z',Z'\big) <\frac{1}{3} \, \gamma .  \medskip \\
\hspace{-0.2cm} {\rm (iii)} & \hspace{-0.25cm}\text{The}\; L\text{-derivatives
of }\ell \text{\ are bounded by } \\
&\; C \, (1+\left\vert y\right\vert _{H}+\left\vert Y\right\vert
_{H}+\left\Vert z\right\Vert _{L_{2}\left( E_{2},H\right) }+\left\Vert
Z\right\Vert _{L_{2}\left( E_{1},H\right) }+\bar{w}_{2}\left( \mu
,\delta _{0}\right) ). \medskip \\
\hspace{-0.15cm} {\rm (iv)} & \hspace{-0.25cm}\text{The}\; L 
\text{-derivatives of}\;\varphi \text{ and } \psi \;\text{are bounded by} \; C\, (1+\left\vert y\right\vert _{H}+\bar{w}_{2}\left( \upsilon ,\delta
_{0}\right) ) \\
& \hspace{-0.2cm} \text{and}\;C \, (1+\left\vert Y\right\vert _{H}+
\bar{w}_{2}\left( \upsilon ,\delta _{0}\right)), \text{
respectively}.
\end{array}
\right.
\end{equation*}

$\left(\mathbf{A}_{6}\right)$:
Denoting $A\left(  t,v,\mu \right)  =\left(-F,f,-G,g\right)  \left(  t,v,\mu \right)$ as  in hypothesis $\left(\mathbf{A}_{2}\right)$, we assume that either $c>0$ and $A$ satisfies $\left(\mathbf{A}_{2}\right)$, or $c<0$ with $A$ satisfying
 $\left(\mathbf{A}_{2}\right)'$.

\bigskip

\smallskip

As we saw in Section~\ref{sec:3}, the condition $0 < \gamma < \frac{1}{2}$ crucial to guarantee the existence of solutions to the adjoint equations of MV-FBDSDE~(\ref{eq:4.2}) in Theorem~\ref{thm:existence for adjoint eqn} below, which is one of the main theorems of this subsection. For additional clarification and similar discussions, refer to Remark~4.2~(i) in \cite{Al-G-relaxed}.

The following theorem addresses the existence and uniqueness of the
solution of MV-FBDSDEs~(\ref{eq:4.2}).
\begin{theorem}\label{thm:4.1} 
For any given admissible control $u_{\cdot}$, if
assumptions $\left(\mathbf{A}_{4}\right)$--$\left(\mathbf{A}_{6}\right)$ hold, then system~(\ref{eq:4.2}) possesses a unique solution.
\end{theorem}
\begin{proof}
Considering that $C^1$ mappings with bounded derivatives are globally Lipschitz, it is evident and straightforward to verify that assumptions $\left(\mathbf{A}_{4}\right)$--$\left(\mathbf{A}_{6}\right)$ imply 
$\left( \mathbf{A}_{1}\right)$--$\left(\mathbf{A}_{3}\right)$. For instance,  
$\left(\mathbf{A}_{4}\right)$(i, iii) and  
$\left(\mathbf{A}_{5}\right)$(ii), along with the definitions of $\bar{w}$ in 
(\ref{eq:2.2a}) and the lifted functions, imply 
 $\left(\mathbf{A}_{1}\right)$(ii, iii). As a result, the proof of the theorem can be derived from Theorem~\ref{thm:main thm1} and 
Theorem~\ref{thm:main thm2}.
\end{proof}

\medskip

Let $u_{\cdot}$ be an arbitrary element of $\mathcal{U}_{ad}$, and let
 $(y_{t}^{u_{\cdot}},Y_{t}^{u_{\cdot}},z_{t}^{u_{\cdot}},Z_{t}^{u_{\cdot}})
$ be the corresponding solution of system (\ref{eq:4.2}). Suppose that 
$\left( \mathbf{A}_{4}\right) $--$\left(\mathbf{A}_{6}\right)$ hold. 
First, we want to introduce the adjoint equations of the MV-FBDSDEs~(\ref{eq:4.2}), and then we present our main result regarding the maximum principle for the optimal control of system~(\ref{eq:4.2}). To this end, let us
define the \emph{Hamiltonian}:
$$\mathcal{H}: [0,T]\times \Omega \times \mathbb{H}^{2}\times K 
\times \mathbb{H}^{2}\times \mathcal{P}_{2}\left( \mathbb{H}^{2}\right) \to \mathbb{R}$$ by the formula:
\begin{equation*}
\begin{array}{ll}
\mathcal{H}(t,y,Y,z,Z,v,p,P,q,Q,\mu ):=\left\langle p,F(t,y,Y,z,Z,v,\mu )\right\rangle
&  \\
\hspace{4.5cm}-\left\langle P,f(t,y,Y,z,Z,v,\mu )\right\rangle +\left\langle
q,G(t,y,Y,z,Z,v,\mu )\right\rangle  &  \\
\hspace{4.5cm}-\left\langle Q,g(t,y,Y,z,Z,v,\mu )\right\rangle -\ell
(t,y,Y,z,Z,v,\mu ). &
\end{array}
\end{equation*}

Using the notation in Section~\ref{sec:appl}, the adjoint equations 
of MV-FBDSDEs~(\ref{eq:4.2}) are the following MV-FBDSDEs:
\begin{equation}\label{eq:4.4}
\left\{
\begin{array}{ll}
dp_{t}^{u_{\cdot}}=\nabla_{Y} \mathcal{H} \big( t,V_{t}^{u_{\cdot}},
u_{t },\chi_{t}^{u_{\cdot}},\mathbb{P}_{V_{t}^{u_{\cdot}}}\big) +\mathbb{\tilde{E}}
\left[ \partial_{\mu_{_Y}}\!
\mathcal{H}\big( t,\tilde{V}_{t}^{u_{\cdot}},u_{t },
\tilde{\chi}_{t}^{u_{\cdot}},\mathbb{P}_{V_{t}^{u_{\cdot}}}\big) \big( V_{t}^{u_{\cdot}}\big)
\right] \Big) dt &  \medskip \\
\hspace{0.75cm} + \left( \nabla_{Z} \mathcal{H}\big( t,V_{t}^{u_{\cdot}},u_{t },\chi_{t}^{u_{\cdot}},\mathbb{P}_{V_{t}^{u_{\cdot}}}\big)
 + \mathbb{\tilde{E}}\left[ \partial_{\mu_{_Z}}\! \mathcal{H}
 \big(t,\tilde{V}_{t}^{u_{\cdot}},u_{t },\tilde{\chi}_{t}^{u_{\cdot}},
 \mathbb{P}_{V_{t}^{u_{\cdot}}}\big) \big( V_{t}^{u_{\cdot}}\big) \right] \right)
dW_{t}  
& \smallskip \\
\hspace{0.75cm} -\, q_{t}^{u_{\cdot}}\, \overleftarrow{dB}_{t}, 
& \medskip \\
dP_{t}^{u_{\cdot}}=\Big(\nabla_{y} \mathcal{H} \big( t,V_{t}^{u_{\cdot}},
u_{t },\chi_{t}^{u_{\cdot}},\mathbb{P}_{V_{t}^{u_{\cdot}}}\big) +\mathbb{\tilde{E}}
\left[ \partial_{\mu_{_y}}\!
\mathcal{H}\big( t,\tilde{V}_{t}^{u_{\cdot}},u_{t },
\tilde{\chi}_{t}^{u_{\cdot}},\mathbb{P}_{V_{t}^{u_{\cdot}}}\big) \big( V_{t}^{u_{\cdot}}\big)\right] \Big) dt &  \smallskip \\
\hspace{0.75cm} + \left( \nabla_{z} \mathcal{H} \big( t,V_{t}^{u_{\cdot}},
u_{t },\chi_{t}^{u_{\cdot}},\mathbb{P}_{V_{t}^{u_{\cdot}}}\big)
 + \mathbb{\tilde{E}}\left[ \partial_{\mu_{_z}}\!
 \mathcal{H}\big( t,\tilde{V}_{t}^{u_{\cdot}},u_{t },\tilde{\chi}_{t}^{u_{\cdot}},
 \mathbb{P}_{V_{t}^{u_{\cdot}}}\big) \big( V_{t}^{u_{\cdot}}\big) \right] \right)
\overleftarrow{dB}_{t} 
&  \smallskip \\
\hspace{0.75cm}
+ \, Q_{t}^{u_{\cdot}} \, dW_{t}, \quad t \in [0,T], &  \medskip \\
p_{0}^{u_{\cdot}}
=- \nabla_{Y} \psi \big( Y_{0}^{u_{\cdot}},\mathbb{P}_{Y_{0}^{u_{\cdot}}}\big) 
-\mathbb{\tilde{E}}\left[ 
\partial_{\mu_{_Y}} \psi \big( \tilde{Y}_{0}^{u_{\cdot}},
\mathbb{P}_{Y_{0}^{u_{\cdot}}}\big) \left( Y_{0}^{u_{\cdot}}\right) \right] ,
& \medskip \\
P_{T}^{u_{\cdot}}= \nabla_{y} \varphi \left(
y_{T}^{u_{\cdot}},\mathbb{P}_{y_{T}^{u_{\cdot}}}\right) +\mathbb{\tilde{E}}
\left[ \partial_{\mu_{_y}} \varphi \left( \tilde{y}_{T}^{u_{\cdot}},\mathbb{P}_{y_{T}^{u_{\cdot}}}\right) \left(
y_{T}^{u_{\cdot}}\right) \right] - c \, p_{T}^{u_{\cdot}}, &
\end{array}
\right.   
\end{equation}
where $V_{t}^{u_{\cdot}}\triangleq (y_{t}^{u_{\cdot}},Y_{t}^{u_{\cdot}},z_{t}^{u_{\cdot
}},Z_{t}^{u_{\cdot}}),$ $\tilde{V}_{t}^{u_{\cdot}}\triangleq(\tilde{y}_{t}^{u_{\cdot}},\tilde{Y}_{t}^{u_{\cdot}},\tilde{z}_{t}^{u_{\cdot}},\tilde{Z}_{t}^{u_{\cdot}}),$ 
$\chi _{t}^{u_{\cdot}}\triangleq(p_{t}^{u_{\cdot}},P_{t}^{u_{\cdot}},q_{t}^{u_{\cdot}},Q_{t}^{u_{\cdot}})$, and 
$\tilde{\chi}_{t}^{u_{\cdot}}\triangleq(\tilde{p}_{t}^{u_{\cdot}},\tilde{P}_{t}^{u_{\cdot}},\tilde{q}_{t}^{u_{\cdot}},\tilde{Q}_{t}^{u_{\cdot}})$. 

\medskip

Here, $\nabla_{Y} \mathcal{H} ( t,V_{t}^{u_{\cdot}},
u_{t },\chi_{t}^{u_{\cdot}},\mathbb{P}_{V_{t}^{u_{\cdot}}})$ is the gradient, defined using the G\^{a}teaux differential: 
$D\mathcal{H} ( t,Y_{t}^{u_{\cdot}}) (h) = \left \langle \nabla_{Y} \mathcal{H} ( t,Y_{t}^{u_{\cdot}}) \, , h \right \rangle_{H}$ at the point $Y_{t}^{u_{\cdot}}$ in the direction $h\in H$, where
$\mathcal{H} ( t,Y_{t}^{u_{\cdot}}):=\mathcal{H} ( t,y_{t}^{u_{\cdot}},Y_{t}^{u_{\cdot}},z_{t}^{u_{\cdot
}},Z_{t}^{u_{\cdot}}, u_{t },\chi_{t}^{u_{\cdot}},\mathbb{P}_{(y_{t}^{u_{\cdot}},Y_{t}^{u_{\cdot}},z_{t}^{u_{\cdot}},Z_{t}^{u_{\cdot}})})$, etc.

\bigskip

In view of our results in Section~\ref{sec:3}, we observe the following theorem.
\begin{theorem}\label{thm:existence for adjoint eqn} 
Under $\left(\mathbf{A}_{4}\right)$--$\left( \mathbf{A}_{6}\right)$, there exists a unique solution $(p^{u_{\cdot}},P^{u_{\cdot}},q^{u_{\cdot}},Q^{u_{\cdot}})$ of the adjoint equations (\ref{eq:4.4}).
\end{theorem}
\begin{proof}
This system~(\ref{eq:4.4}) can be expressed as a linear system of MV-FBDSDEs on the arguments $(p^{u_{\cdot}},P^{u_{\cdot}},q^{u_{\cdot}},Q^{u_{\cdot}})$. With the assumptions
$\left(\mathbf{A}_{4}\right)$--$\left( \mathbf{A}_{6}\right)$, it is evident that this linear system satisfies $\left(\mathbf{A}{1}\right)$,
 $\left(\mathbf{A}{2}\right)'$, and $\left( \mathbf{A}_{3}\right)$. For more detailed information on a similar approach for FBDSDEs, one can refer to the methodology employed in our previous work \cite{Al-G-relaxed}. Therefore, the desired result follows 
 from Theorem~\ref{thm:main thm1}, Theorem~\ref{thm:main thm2}, 
 and Remark~\ref{rk:rk-on-assumptions}~(ii).
\end{proof}

\bigskip

Now, we present the main theorem of this section.
\begin{theorem}[Sufficient conditions for optimality] \label{theorem 4.1}
Assume that conditions $\left(\mathbf{A}_{4}\right)$--$\left(\mathbf{A}_{6}\right)$ hold. Given 
$\widehat{u}_{\cdot }\in \mathcal{U}_{ad},$ let 
$V_{t}^{\widehat{u}_{\cdot}}\equiv(y_{t}^{\widehat{u}_{\cdot }},Y_{t}^{
\widehat{u}_{\cdot }},z_{t}^{\widehat{u}_{\cdot }},Z_{t}^{\widehat{u}_{\cdot }})$ and
$\left( p^{\widehat{u}},P^{\widehat{u}},q^{\widehat{u}},Q^{\widehat{u}}\right)$ be the corresponding solutions of MV-FBDSDEs~(\ref{eq:4.2}) and (\ref{eq:4.4}), respectively. Assume the following:
\medskip \\ 
{\rm (i)}  $\varphi $ and $\psi $ are convex. 
\medskip \\
{\rm (ii)} For all $t \in [0,T],\;\mathbb{P}-a.s.$, the function $\mathcal{H}(t,\cdot ,\cdot ,\cdot ,\cdot ,\cdot ,\cdot ,p^{\widehat{u}},P^{\widehat{u}},q^{\widehat{u}},Q^{\widehat{u}})$ is
concave. \medskip \\
{\rm (iii)}  We have
\begin{equation*}\label{eq:4.5}
\mathcal{H}\big(t,V_{t}^{\widehat{u}_{\cdot}},\widehat{u}_{t},\chi_{t}^{\widehat{u}_{\cdot }},\mathbb{P}_{V_{t}^{\widehat{u}_{\cdot }}}\big)=\max_{v \in K}\mathcal{H}\big(t,V_{t}^{\widehat{u}_{\cdot
}},v,\chi _{t}^{\widehat{u}_{\cdot }},\mathbb{P}
_{V_{t}^{\widehat{u}_{\cdot }}}\big),\;  \; \, a.e.\; t, \;\mathbb{P}-a.s.  
\end{equation*}
Then
 $(y^{\widehat{u}_{\cdot}},Y^{\widehat{u}_{\cdot}},z^{\widehat{u}_{\cdot}},
 Z^{\widehat{u}_{\cdot}},\widehat{u}_{\cdot})$
is an optimal solution of the control problem (\ref{eq:4.1})-(\ref{eq:4.3}).
\end{theorem}
\begin{proof}
Let $\widehat{u}_{\cdot }\in \mathcal{U}_{ad}$  be an arbitrary candidate for an optimal control. Consider its associated trajectory 
$V_{t}^{\widehat{u}_{\cdot}}\triangleq (y_{t}^{\widehat{u}_{\cdot }},Y_{t}^{\widehat{u}_{\cdot }},z_{t}^{\widehat{u}_{\cdot }},Z_{t}^{\widehat{u}_{\cdot }})$. For any $u_{\cdot}\in \mathcal{U}_{ad}$ with its associated trajectory 
$V_{t}^{u_{\cdot}}\triangleq (y_{t}^{u_{\cdot}},Y_{t}^{u_{\cdot}},z_{t}^{u_{\cdot}},Z_{t}^{u_{\cdot}})$, we have 
\begin{align*}
& J(u_{\cdot})-J(\hat{u}_{\cdot})=\mathbb{E}\left[  \varphi(y_{T}^{u_{\cdot}},\mathbb{P}_{y_{T}^{u}})-\varphi(y_{T}^{\hat{u}_{\cdot}},\mathbb{P}_{y_{T}^{\widehat{u}_{\cdot}}})\right]  +\mathbb{E}\left[  \psi(Y_{0}^{u_{\cdot}},\mathbb{P}_{Y_{0}^{u_{\cdot}}})-\psi(Y_{0}^{\widehat{u}_{\cdot}},\mathbb{P}_{Y_{0}^{\widehat{u}_{\cdot}}})\right]  \\
& \hspace{3cm}+\, \mathbb{E}\left[  \int_{0}^{T}\left(  \ell(t,V_{t}^{u_{\cdot}},u_{t},\mathbb{P}_{ V_{t}^{u_{\cdot}}})-\ell(t,V_{t}^{\widehat{u}_{\cdot}},\widehat{u}_{t},\mathbb{P}_{V_{t}^{\widehat{u}_{\cdot}} })\right)  dt\right] .
\end{align*}

In accordance with the convexity of $\varphi$ and $\psi$ as indicated in (\ref{eq:convexity}), we deduce 
\begin{align*}
& \mathbb{E}\left[  \varphi(y_{T}^{u_{\cdot}},\mathbb{P}_{y_{T}^{u}}
)-\varphi(y_{T}^{\hat{u}_{\cdot}},\mathbb{P}_{y_{T}^{\widehat{u}_{\cdot}}
})\right]  \geq \mathbb{E}\left[   \big\langle \nabla_{y} \, \varphi 
\big(y_{T}^{\widehat{u}_{\cdot}},\mathbb{P}_{y_{T}^{\widehat{u}_{\cdot}}}\big)
,y_{T}^{u_{\cdot}}-y_{T}^{\widehat{u}_{\cdot}} \big\rangle_{H} \, \right]
\\
& \hspace{1cm}
+  \mathbb{E}\left[ \mathbb{\tilde{E}}\big[  \big\langle
  \partial_{\mu_{y}} \varphi \big(
y_{T}^{\widehat{u}_{\cdot}},\mathbb{P}_{y_{T}^{\widehat{u}_{\cdot}}}\big)
\big(  \tilde{y}_{T}^{\widehat{u}_{\cdot}}\big)  ,\tilde{y}
_{T}^{u_{\cdot}}-\tilde{y}_{T}^{\widehat{u}_{\cdot}} \big\rangle_{H} \, \big] \right]
\\
& \hspace{1cm} =\mathbb{E}\left[  \left \langle \nabla_{y} \, \varphi \big(
y_{T}^{\widehat{u}_{\cdot}},\mathbb{P}_{y_{T}^{\widehat{u}_{\cdot}}}\big)
+\mathbb{\tilde{E}}\left[ \partial_{\mu_{y}} \varphi \big(
\tilde{y}_{T}^{\widehat{u}_{\cdot}},\mathbb{P}_{y_{T}^{\widehat{u}_{\cdot}}
}\big)  \big(  y_{T}^{\widehat{u}_{\cdot}}\big)  \right]
  ,y_{T}^{u_{\cdot}}-y_{T}^{\widehat{u}_{\cdot}}\right \rangle_{H} \, \right]  ,
\end{align*}
and
\begin{align*}
& \mathbb{E}\left[  \psi(Y_{0}^{u_{\cdot}},\mathbb{P}_{Y_{0}^{u_{\cdot}}
})-\psi(Y_{0}^{\widehat{u}_{\cdot}},\mathbb{P}_{Y_{0}^{\widehat{u}_{\cdot}}
})\right]  \geq \mathbb{E}\left[  \big\langle \nabla_{_Y} \psi \big(
Y_{0}^{\widehat{u}_{\cdot}},\mathbb{P}_{Y_{0}^{\widehat{u}_{\cdot}}}\big)
,Y_{0}^{u_{\cdot}}-Y_{0}^{\widehat{u}_{\cdot}} \big\rangle_{H} \, \right]
\\
& \hspace{1cm}
+ \mathbb{E}\left[  \mathbb{\tilde{E}} \big[ \big\langle
  \partial_{\mu_{_Y}} \psi \big(
Y_{0}^{\widehat{u}_{\cdot}},\mathbb{P}_{Y_{0}^{\widehat{u}_{\cdot}}}\big)
\big(  \tilde{Y}_{0}^{v_{\cdot}}\big)  ,\tilde{Y}_{0}^{u_{\cdot}
}-\tilde{Y}_{0}^{\widehat{u}_{\cdot}} \big\rangle_{H}  \, \big] \right]  
\\
& \hspace{1cm} =\mathbb{E}\left[  \left \langle \nabla_{_Y} \psi \big(
Y_{0}^{\widehat{u}_{\cdot}},\mathbb{P}_{Y_{0}^{\widehat{u}_{\cdot}}}\big)
+ \mathbb{\tilde{E}} \Big[ \partial_{\mu_{_Y}} \psi \big(
\tilde{Y}_{0}^{\widehat{u}_{\cdot}},\mathbb{P}_{Y_{0}^{\widehat{u}
_{\cdot}}}\big)  \left(  Y_{0}^{v_{\cdot}}\right)  \Big]  ,Y_{0}
^{u_{\cdot}}-Y_{0}^{\widehat{u}_{\cdot}}\right \rangle_{H} \, \right] .
\end{align*}
Therefore, by implementing (\ref{eq:4.4}), it follows that
\begin{eqnarray}\label{eq:4.6}
&&\hspace{-0.75cm}J\left( u_{\cdot
}\right) -J\left(\widehat{u}_{\cdot }\right) \geq \mathbb{E} \, 
\big[\langle P_{T}^{\widehat{u}_{\cdot}} + c \, p_{T}^{\widehat{u}_{\cdot}} ,
y_{T}^{u_{\cdot}}-y_{T}^{\widehat{u}_{\cdot}}\rangle_{H} \big]-\mathbb{E}\, 
\big[\langle p_{0}^{\widehat{u}_{\cdot}},Y_{0}^{u_{\cdot}}-Y_{0}^{\widehat{u}_{\cdot }}\rangle_{H} \big]  \nonumber \\
&&\hspace{2.2cm}+ \, \mathbb{E}\left[\int_{0}^{T}\Big(\ell (t,V_{t}^{u_{\cdot}},u_{t},\mathbb{P}_{V_{t}^{u_{\cdot}}})-\ell(t,V_{t}^{\widehat{u}_{\cdot}},\widehat{u}_{t},\mathbb{P}_{V_{t}^{\widehat{u}_{\cdot}}})\Big)dt\right]. 
\end{eqnarray}

Next, by calculating $\mathbb{E}\,\big[\langle P_{T}^{\widehat{u}_{\cdot }},y_{T}^{u_{\cdot}}-y_{T}^{\widehat{u}_{\cdot}}\rangle_{H} \big]$ and
 $\mathbb{E}\,\big[\langle p_{0}^{\widehat{u}_{\cdot}},Y_{0}^{u_{\cdot}}-Y_{0}^{\widehat{u}_{\cdot }}\rangle_{H} \big]$ through applying It\^{o}'s formula (e.g. as in \cite{CD2} or \cite{Gozzi}) to compute 
 $\langle P_{t}^{\widehat{u}_{\cdot}},y_{t}^{u_{\cdot}}-y_{t}^{\widehat{u}_{\cdot}}\rangle_{H}$ and 
 $\langle p_{t}^{\widehat{u}_{\cdot}},Y_{t}^{u_{\cdot}}-Y_{t}^{\widehat{u}_{\cdot}}\rangle_{H}$,
respectively (cf. e.g. \cite{Al-G-relaxed}), and then utilizing the equality 
$$\mathbb{E} \, \big[\langle c \, p_{T}^{\widehat{u}_{\cdot}} ,
y_{T}^{u_{\cdot}}-y_{T}^{\widehat{u}_{\cdot}}\rangle_{H} \big] = \mathbb{E}\,\big[\langle p_{T}^{\widehat{u}_{\cdot }},Y_{T}^{u_{\cdot}}-Y_{T}^{\widehat{u}_{\cdot}}\rangle_{H} \big] ,$$
which is evident from (\ref{eq:4.2}), we ultimately obtain 
\begin{equation*}
\begin{array}{ll}
J(u_{\cdot})-J(\widehat{u}_{\cdot }) 
\geq \mathbb{E}\left[\displaystyle\int_{0}^{T}\left\langle \nabla_{y} \mathcal{H}\big(t,V_{t}^{\widehat{u}_{\cdot}},\widehat{u}_{t},\chi_{t}^{\widehat{u}_{\cdot }},\mathbb{P}_{V_{t}^{\widehat{u}_{\cdot }}}\big) ,y_{t}^{u_{\cdot
}}-y_{t}^{\widehat{u}_{\cdot }}\right\rangle_{H} dt\right] &   \medskip \\
\hspace{2cm} + \, \mathbb{E}\left[\displaystyle\int_{0}^{T}\left\langle \mathbb{\tilde{E}}\left[ \partial_{\mu_{_y}}\! \mathcal{H}\big( t,\tilde{V}_{t}^{\widehat{u}_{\cdot}},\widehat{u}_{t},\tilde{\chi}_{t}^{\widehat{u}_{\cdot }},\mathbb{P}
_{V_{t}^{\widehat{u}_{\cdot }}}\big) \left( V^{\widehat{u}_{\cdot }}\right) \right] ,y_{t}^{u_{\cdot
}}-y_{t}^{\widehat{u}_{\cdot }}\right\rangle_{H} dt\right] &   \medskip \\
\hspace{2cm}+ \,\mathbb{E}\left[\displaystyle\int_{0}^{T}\left\langle \nabla_{Y}\mathcal{H}\big(t,V_{t}^{\widehat{u}_{\cdot}},\widehat{u}_{t},\chi_{t}^{\widehat{u}_{\cdot }},\mathbb{P}_{V_{t}^{\widehat{u}_{\cdot }}}\big) ,Y_{t}^{u_{\cdot
}}-Y_{t}^{\widehat{u}_{\cdot }}\right\rangle_{H} dt\right] &   \medskip \\
\hspace{2cm}+ \,\mathbb{E}\left[\displaystyle\int_{0}^{T}\left\langle \mathbb{\tilde{E}}\left[ \partial_{\mu_{_Y}}\! \mathcal{H}\big( t,\tilde{V}_{t}^{\widehat{u}_{\cdot}},\widehat{u}_{t},\tilde{\chi}_{t}^{\widehat{u}_{\cdot }},\mathbb{P}_{V_{t}^{\widehat{u}_{\cdot }}}\big) \big( V^{\widehat{u}_{\cdot }}\big) \right] ,Y_{t}^{u_{\cdot
}}-Y_{t}^{\widehat{u}_{\cdot }}\right\rangle_{H} dt\right] &   \medskip \\
\hspace{2cm}+ \,\mathbb{E}\left[\displaystyle\int_{0}^{T}\left\langle \nabla_{z} \mathcal{H}\big(t,V_{t}^{\widehat{u}_{\cdot}},\widehat{u}_{t},\chi_{t}^{\widehat{u}_{\cdot }},\mathbb{P}_{V_{t}^{\widehat{u}_{\cdot }}}\big) ,z_{t}^{u_{\cdot
}}-z_{t}^{\widehat{u}_{\cdot }}\right\rangle_{H} dt\right] &   \medskip \\
\hspace{2cm}+ \,\mathbb{E}\left[\displaystyle\int_{0}^{T}\left\langle \mathbb{\tilde{E}}\left[ \partial_{\mu_{_z}}\! \mathcal{H}\big( t,\tilde{V}_{t}^{\widehat{u}_{\cdot}},\widehat{u}_{t},\tilde{\chi}_{t}^{\widehat{u}_{\cdot }},\mathbb{P}
_{V_{t}^{\widehat{u}_{\cdot }}}\big) \big( V^{\widehat{u}_{\cdot }}\big) \right] ,z_{t}^{u_{\cdot
}}-z_{t}^{\widehat{u}_{\cdot }}\right\rangle_{H} dt\right] &   \medskip \\
\hspace{2cm}+ \,\mathbb{E}\left[\displaystyle\int_{0}^{T}\left\langle \nabla_{Z} \mathcal{H}\big(t,V_{t}^{\widehat{u}_{\cdot}},\widehat{u}_{t},\chi_{t}^{\widehat{u}_{\cdot }},\mathbb{P}_{V_{t}^{\widehat{u}_{\cdot }}}\big) ,Z_{t}^{u_{\cdot
}}-Z_{t}^{\widehat{u}_{\cdot }}\right\rangle_{H} dt\right] &  \medskip \\
\hspace{2cm}+ \,\mathbb{E}\left[\displaystyle\int_{0}^{T}\left\langle \mathbb{\tilde{E}}\left[ \partial_{\mu_{_Z}}\! \mathcal{H}\big( t,\tilde{V}_{t}^{\widehat{u}_{\cdot}},\widehat{u}_{t},\tilde{\chi}_{t}^{\widehat{u}_{\cdot }},\mathbb{P}
_{V_{t}^{\widehat{u}_{\cdot }}}\big) \big( V^{\widehat{u}_{\cdot }}\big) \right] ,Z_{t}^{u_{\cdot
}}-Z_{t}^{\widehat{u}_{\cdot }}\right\rangle_{H} dt\right]
&   \medskip \\
\hspace{2cm} - \, \mathbb{E}\left[\displaystyle\int_{0}^{T}\Big(\mathcal{H}\big(t,V
_{t}^{u_{\cdot}},u_{t},\chi_{t}^{u_{\cdot}},\mathbb{P}_{V_{t}^{u_{\cdot}}}\big)-\mathcal{H}\big(t,V_{t}^{\widehat{u}_{\cdot}},\widehat{u}_{t},\chi_{t}^{\widehat{u}_{\cdot}},\mathbb{P}
_{V_{t}^{\widehat{u}_{\cdot }}}\big)\Big)dt\right].
&
\end{array}
\end{equation*}

Using the concavity of $\mathcal{H}$ assumed in (ii), it follows in a standard way that
\begin{equation}\label{eq:4.7}
J(u_{\cdot})-J(\widehat{u}_{\cdot })\geq - \, \mathbb{E}\left[\int_{0}^{T}\left\langle 
\nabla_{u} \mathcal{H}\big(t,V_{t}^{\widehat{u}_{\cdot}},\widehat{u}_{t},\chi_{t}^{\widehat{u}_{\cdot }},\mathbb{P}
_{V_{t}^{\widehat{u}_{\cdot }}}\big) ,u_{t}-\widehat{u}_{t }\right\rangle_{K} dt\right]. 
\end{equation}
However, the condition (iii) implies that the function $v\mapsto \mathcal{H}\big(t,V_{t}^{\widehat{u}_{\cdot
}},v,\chi _{t}^{\widehat{u}_{\cdot }},\mathbb{P}
_{V_{t}^{\widehat{u}_{\cdot }}}\big)$ is maximal for $v=\widehat{u}_{t}$, so we have
\begin{equation*}
\begin{array}{ll}
\left\langle \nabla_{u} \mathcal{H}\big(t,V_{t}^{\widehat{u}_{\cdot}},\widehat{u}_{t},\chi_{t}^{\widehat{u}_{\cdot }},\mathbb{P}
_{V_{t}^{\widehat{u}_{\cdot }}}\big) ,u_{t}-\widehat{u}_{t}\right\rangle_{K} \leq 0, \quad a.e.\; t,\; \mathbb{P}-a.s. &
\end{array}
\end{equation*}
Consequently, (\ref{eq:4.7}) simplifies to
\begin{equation*}
\begin{array}{ll}
J(u_{\cdot})-J(\widehat{u}_{\cdot })\geq 0 . &
\end{array}
\end{equation*}
Since $u_{\cdot}$ is an arbitrary admissible control, we conclude that 
$\widehat{u}_{\cdot}$ is an optimal control to the control problem 
(\ref{eq:4.1})--(\ref{eq:4.3}).
\end{proof}

\end{document}